\pdfoutput=1
\documentclass[12pt]{article}
\usepackage[margin=1in]{geometry}
\usepackage{amsmath,amssymb,amsthm}
\usepackage{hyperref}
\usepackage{booktabs}
\usepackage{longtable}
\usepackage{array}
\usepackage{enumitem}
\usepackage{tikz}
\usepackage{graphicx}
\usepackage{float}
\usepackage{caption}
\usepackage{xcolor}
\usepackage{mathpazo}
\usepackage{fancyhdr}
\fancypagestyle{plain}{\fancyhf{}\fancyfoot[C]{\small \textcopyright\ 2026 The Mereon Legacy CIC. All rights reserved.}\fancyfoot[R]{\thepage}}
\usetikzlibrary{patterns}

\theoremstyle{definition}

\title{\textbf{The Mereon System, the 600-Cell, and the\\
Exceptional Algebras $E_6$, $E_7$, $E_8$:\\
Exact Correspondence via $H_3 \subset H_4$ Symmetry\\
and the Eigenform Loop}}

\author{
Robert W.\ Gray\thanks{Independent researcher.} \and
Lynnclaire Dennis\thanks{The Mereon Legacy CIC.} \and
Louis H.\ Kauffman\thanks{Department of Mathematics, Statistics and Computer Science, University of Illinois at Chicago, 851 South Morgan Street, Chicago, IL 60607-7045, USA; and International Institute for Sustainability with Knotted Chiral Meta Matter (WPI-SKCM$^2$), Hiroshima University, 1-3-1 Kagamiyama, Higashi-Hiroshima, Hiroshima 739-8531, Japan.}
}

\date{Preprint Draft v07 --- 31 March 2026}

\begin{document}
\maketitle

\begin{abstract}
The Mereon System is a geometric framework of nested polyhedra in
$\mathbb{R}^3$: a crystallographic 144-face core (M144p, symmetry $O_h$),
a non-crystallographic 120-face boundary (M120p, symmetry $H_3$), and an
intermediate focusing sphere.

We identify the unit 3-sphere $S^3 \subset \mathbb{H} \cong \mathbb{R}^4$ with the group of
unit quaternions $\mathrm{Sp}(1)$. Every M120p vertex direction in
$\mathbb{R}^3 \cong \mathrm{Im}(\mathbb{H})$ is the stereographic image
of at least one element of the binary icosahedral group
$2I \subset \mathrm{Sp}(1)$~\cite{ChoiLee2018}: all 62 vertex directions
are matched exactly, verified vertex by vertex (all 62 of 62 vertex
directions matched). The stereographic projection of the 600-cell maps all 120
elements of $2I$ onto these 62 directions, organised into 8 type-pure
shells at radii determined by the golden ratio $\varphi$. The M120p is
thus the $H_3$ shadow of the $H_4$-symmetric 600-cell. The three vertex
types correspond to three latitudes on $S^3$, and the vertex classification
reflects a field-extension cascade of the binary group quaternion
coordinates:
$\mathbb{Q} \to \mathbb{Q}(\sqrt{2}) \to \mathbb{Q}(\sqrt{5})$
(where $\mathbb{Q}(\sqrt{2})$ refers to the unit quaternion components of
$2O$, not to the integer coordinates of the M144p).

The 24 elements of $2I$ at the highest non-trivial latitude
($|w| = \varphi/2$) project to an inner icosahedron at radius
$1/\sqrt{4\varphi+3}$---the innermost non-trivial shell of the 600-cell
projection. This inner icosahedron is exactly radially aligned with the 12
C-vertices of the M120p, at a radius ratio of exactly $1/\varphi$. Its
12 vertices coincide with the focusing sphere of the Mereon System,
first described by Dennis~\cite{Dennis1997}.

The three exceptional Lie algebras $E_6$, $E_7$, $E_8$ are realised
geometrically: $E_6 \subset E_8$ via the McKay correspondence
\cite{McKay1980} applied to the inclusion $2T \subset 2I$, visible as
the 24-cell inscribed in the 600-cell; and $E_7$ independently via the
M144p's crystallographic $O_h$ symmetry (McKay applied to $2O$). The
focusing sphere marks the geometric boundary between these two regimes.
The Brieskorn variety $z_1^2 + z_2^3 + z_3^5 = 0$~\cite{Brieskorn1966}
closes the structure into an eigenform loop~\cite{Kauffman2005}
connecting the trefoil knot, $M(2,3,5)$~\cite{MilnorTriangle}, $2I$,
and the M120p.
\end{abstract}

\section{Introduction}
\label{sec:intro}

Our most recent work concerns how the three-dimensional polyhedral Mereon structure
(the 120 polyhedron) is the precise projection from four-space of the 600-cell, an analogue
in four-dimensional space of a regular solid. The 600-cell is made from 120 copies of a dodecahedron that are fitted together so that each dodecahedral face is matched to the face
of another dodecahedron (much as the pentagonal faces of the dodecahedron are matched
along their edges). Thus this essential part of the Mereon structure is a projection from
a higher-dimensional space of an even more symmetrical entity. The theme that three-dimensional structures, earthly structures, networked structures, structures involved in
our understanding and communication, would be or should be seen as projections from a
higher-dimensional whole is part of perennial philosophy. Here we are seeing an instantiation of this theme and the dreams with which it is allied. The 600-cell and its associated
geometries have been studied for some time by mathematicians and by physicists for relations with geometry, topology, knot theory, particle physics and even cosmology and
string theory. It is more than exciting that there is a direct connection of the Mereon System with the 600-cell and the wide-ranging conversation with which it is associated. We
expect much more from this connection as the search goes on.

The present paper makes this connection precise. We establish the exact correspondence between the Mereon System and the 600-cell, and show how the nested architecture realises all three exceptional Lie algebras $E_6$, $E_7$, $E_8$. The mathematical framework is as follows.

The Coxeter groups $H_3$ and $H_4$ \cite{Humphreys1990} are the symmetry groups of the
icosahedron and dodecahedron (in 3 dimensions) and the 600-cell and 120-cell
(in 4 dimensions), respectively. They are the \emph{only}
non-crystallographic Coxeter groups in their respective dimensions---no
lattice in $\mathbb{R}^3$ or $\mathbb{R}^4$ carries icosahedral
symmetry. The inclusion $H_3 \subset H_4$ is a classical fact: the
icosahedral symmetry of a 3-dimensional cross-section of the 600-cell
is inherited from the full 4-dimensional symmetry.

The Mereon System, first described by Dennis \cite{Dennis1997},
developed topologically by Kauffman since 1995
\cite{KauffmanKnots2013, GrayDennisKauffman2018} and geometrically
by Gray since 1998 \cite{Gray2002, GrayDennisKauffman2018}, is a
dynamic geometric framework built from nested polyhedra. The dynamics will not be presented
in this paper, but see References \cite{Dennis1997} and \cite{GrayDennisKauffman2018}.
Its boundary, the M120p, carries
$H_3$ symmetry. Its core, the 144-face polyhedron (M144p), carries the
crystallographic octahedral symmetry $O_h$. Between them lies a
focusing sphere coupling two algebraically disjoint number fields.

Several concepts from group theory and topology are central to this paper.
The \emph{binary icosahedral group} $2I$ is the double cover of the
icosahedral rotation group $I \cong A_5$ inside the unit quaternions
$\mathrm{Sp}(1) \cong S^3 \subset \mathbb{H} \cong \mathbb{R}^4$; it has order 120, and its
elements are exactly the 120 vertices of the 600-cell
\cite{ChoiLee2018, ConwaySmith2003}. The \emph{McKay correspondence}
\cite{McKay1980} is a bijection between finite subgroups of
$\mathrm{SU}(2)$ and the affine simply-laced Dynkin diagrams: it maps the
binary tetrahedral, octahedral, and icosahedral groups $2T$, $2O$, $2I$
to $\hat{E}_6$, $\hat{E}_7$, $\hat{E}_8$ respectively. A
\emph{Brieskorn variety} $V(p,q,r) : z_1^p + z_2^q + z_3^r = 0$
\cite{Brieskorn1966} is an algebraic surface whose intersection with a
small sphere around the origin (its \emph{link}) is a 3-manifold encoding
the singularity type; for $(p,q,r) = (2,3,5)$, the link is the Poincar\'{e}
homology sphere $M(2,3,5)$ \cite{MilnorTriangle}. An \emph{eigenform}
\cite{Kauffman2005} is a fixed point of a recursive process: the system
generates the structure whose properties regenerate the system.

\begin{figure}[H]
\centering
\includegraphics[width=0.75\textwidth]{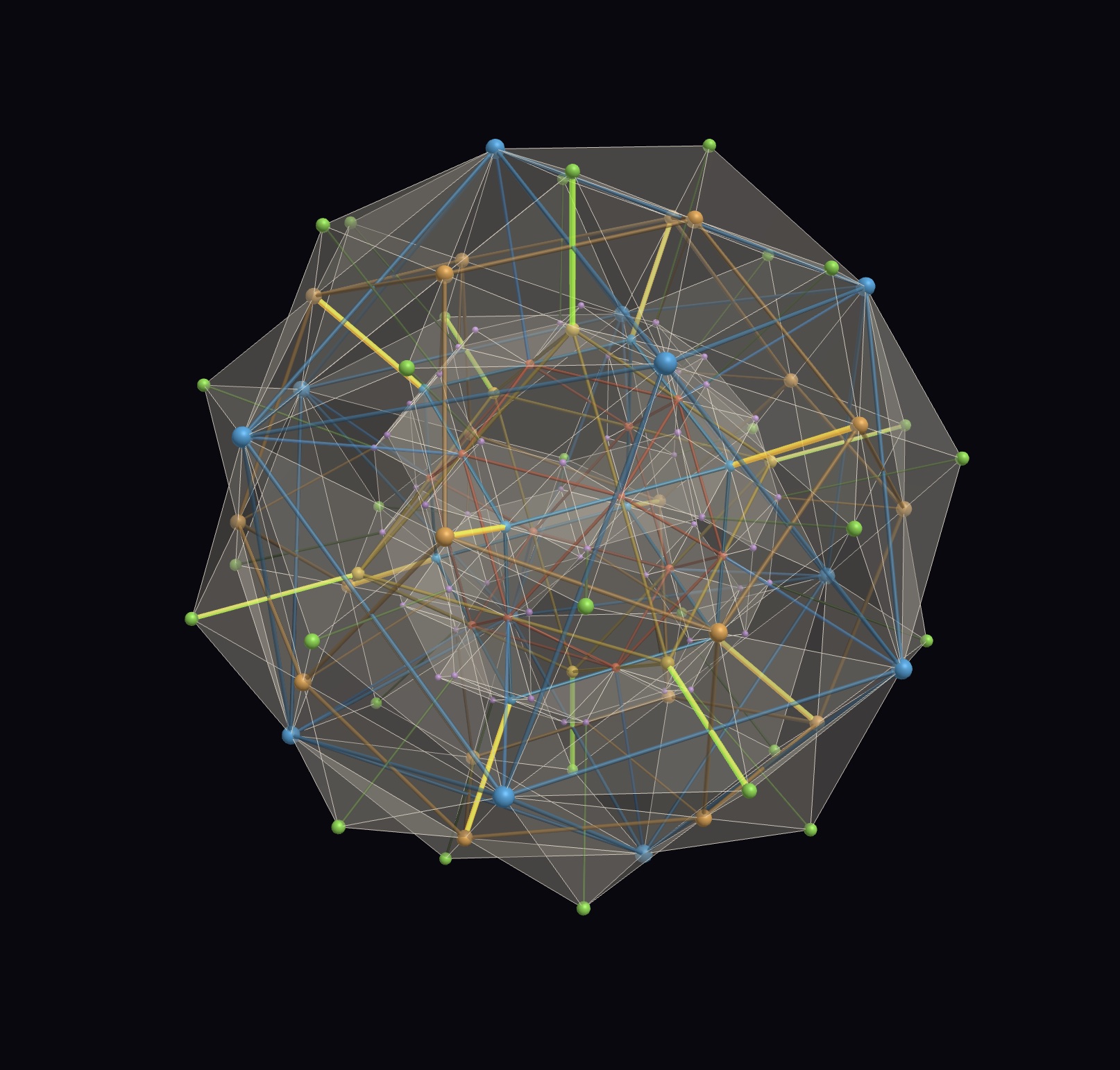}
\caption{The Mereon System: nested geometry with all edges visible.
The M120p boundary (outermost, $H_3$ symmetry) shows its three vertex
types: A-vertices (green, dodecahedral, 3-fold), C-vertices (blue,
icosahedral, 5-fold), and B-vertices (orange, icosidodecahedral,
2-fold). Inside it, the M144p core ($O_h$ symmetry) is visible with
its own edge network and the focusing sphere. The trefoil knot region
lies between the focusing sphere and the M120p boundary (shown
separately in Figure~4).}
\label{fig:nested}
\end{figure}

The overall architecture of the Mereon System is shown schematically
in Figure~\ref{fig:nested}.

This paper establishes that the M120p is the $H_3$ shadow of the
$H_4$-symmetric 600-cell, and that the nested Mereon architecture
realises all three exceptional Lie algebras $E_6$, $E_7$, $E_8$
through the interplay of crystallographic and non-crystallographic
symmetry. We present both the results and the discovery pathway,
including the $w$-coordinate recovery that revealed the fourth dimension
hidden in Gray's three-dimensional coordinates. The connection of the Mereon System
to the ADE classification of Dynkin diagrams is discussed in Section~\ref{sec:discussion}.
The M120p is further shown to be distinct
from the Disdyakis Triacontahedron despite sharing the same 12+20+30
vertex partition.

\subsection{Notation}

$\varphi = (1+\sqrt{5})/2$. Key identities: $\varphi^2 = \varphi + 1$,
$2\varphi - 1 = \sqrt{5}$, $2\varphi + 1 = \varphi^3$,
$4\varphi + 3 = \varphi^2(\varphi^2 + 1)$.
Quaternions: $q = w + xi + yj + zk$.
Stereographic projection from $(-1,0,0,0) \in S^3$:
\begin{equation}\label{eq:stereo}
  \pi(q) = \frac{(x, y, z)}{1+w} \in \mathbb{R}^3.
\end{equation}
Double cover: $R_q(v) = qv\bar{q}$, giving
$\mathrm{SU}(2) \to \mathrm{SO}(3)$, kernel $\{\pm 1\}$.

Vertex naming convention: $\text{A} = \text{Input}$ (20 vertices, 3-fold,
dodecahedral, innermost), $\text{B} = \text{Thruput}$ (30 vertices,
2-fold, icosidodecahedral, outermost), $\text{C} = \text{Output}$
(12 vertices, 5-fold, icosahedral, intermediate).

\section{Methods}
\label{sec:methods}

The mathematical results in this paper were developed collaboratively
by the three authors over a period of months, combining Gray's
long-standing geometric constructions with Kauffman's topological
expertise and Dennis's originating insight into the Mereon System.

Computational exploration was conducted with Claude Opus 4
(Anthropic, claude-opus-4-6) as tool and writing co-operator under author direction.
Mathematical authority rests with the named authors. All algebraic
identities and vertex correspondences were verified independently
by direct computation.

\paragraph{Validation.}
Every claim of exact correspondence (e.g.\ the 62~of~62 vertex match,
the shell radii, the 1/$\varphi$ ratio) was checked by explicit
numerical and symbolic computation. Where stereographic projection
maps are asserted to be exact, both the algebraic derivation and a
floating-point verification to machine precision were performed.
The computational dialogue was validated while cultivating the
questions put to Claude by repeatedly asking for clarification,
elaboration, verification, etc.\ of specific parts of a response.

\section{The Mereon System: Nested Dynamic Polyhedra}
\label{sec:system}

The Mereon System is not a single polyhedron but a system of nested dynamic
structures sharing a common topological language. However, we focus on the static
description of the system in this paper. 

\subsection{The M144p core: crystallographic ($O_h$)}
\label{sec:core}

The innermost structure, referred to here as the M144p, has
74 vertices, 216 edges, and 144 triangular faces. The vertex coordinates can be scaled so that they are all integers. The M144p has octahedral symmetry $O_h$ (order 48).

\begin{figure}[H]
\centering
\begin{minipage}[b]{0.31\textwidth}
  \centering
  \includegraphics[width=\textwidth]{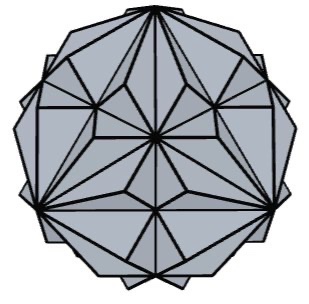}
\end{minipage}
\hfill
\begin{minipage}[b]{0.31\textwidth}
  \centering
  \includegraphics[width=\textwidth]{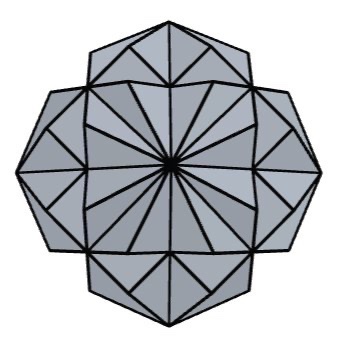}
\end{minipage}
\hfill
\begin{minipage}[b]{0.31\textwidth}
  \centering
  \includegraphics[width=\textwidth]{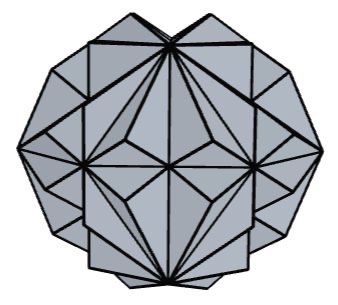}
\end{minipage}
\caption{The M144p 144-face core polyhedron (three views).}
\label{fig:M144p}
\end{figure}

The 74 vertices occupy four Face Centred Cubic lattice (FCC) shells:

\begin{center}
\begin{tabular}{llll}
\toprule
$r^2$ & Coordinate type & Count & Note \\
\midrule
$8$ & $(\pm 2, \pm 2, 0)$ perms & 12 & Edge midpoints \\
$14$ & $(\pm 1, \pm 2, \pm 3)$ perms & 48 & Hexagon vertices \\
$12$ & $(\pm 2, \pm 2, \pm 2)$ & 8 & Face center vertices \\
$16$ & $(\pm 4, 0, 0)$ perms & 6 & 4-f Octa vertices \\
\bottomrule
\end{tabular}
\end{center}

The M144p core is a crystallographic structure compatible with lattice symmetry.

\subsubsection*{Construction of the M144p}

The M144p is constructed from the FCC lattice via the following
sequence of five steps applied to a 6-frequency octahedron (where
\emph{$n$-frequency} means that $n$ steps are required along an edge
to travel from one vertex of the octahedron to an adjacent vertex).

\begin{figure}[H]
\centering
\begin{minipage}[b]{0.31\textwidth}
  \centering
  \includegraphics[width=\textwidth]{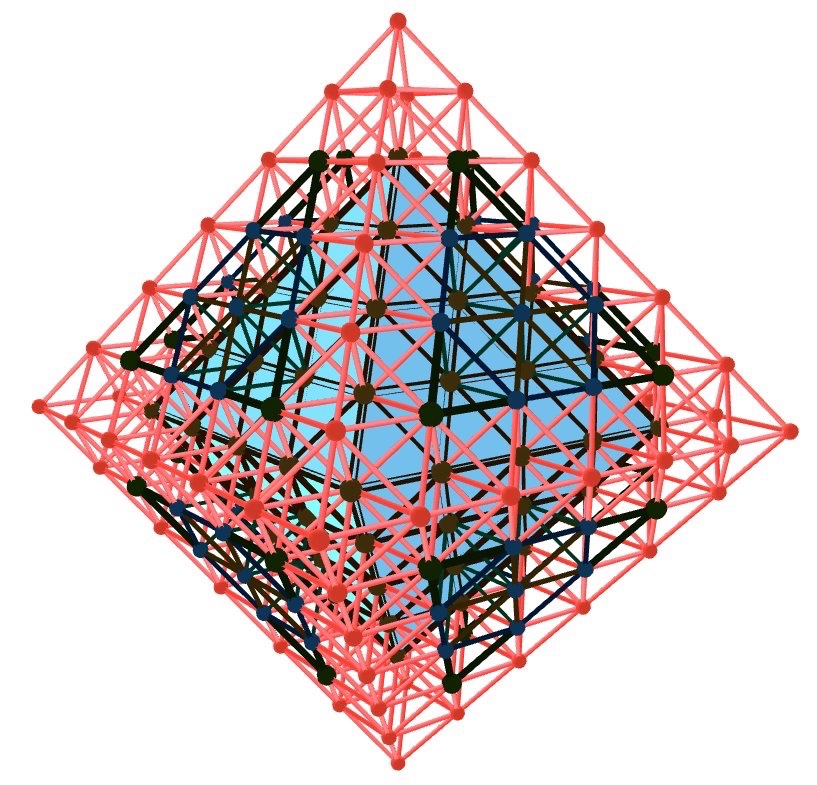}
  \caption*{Step 1}
\end{minipage}
\hfill
\begin{minipage}[b]{0.31\textwidth}
  \centering
  \includegraphics[width=\textwidth]{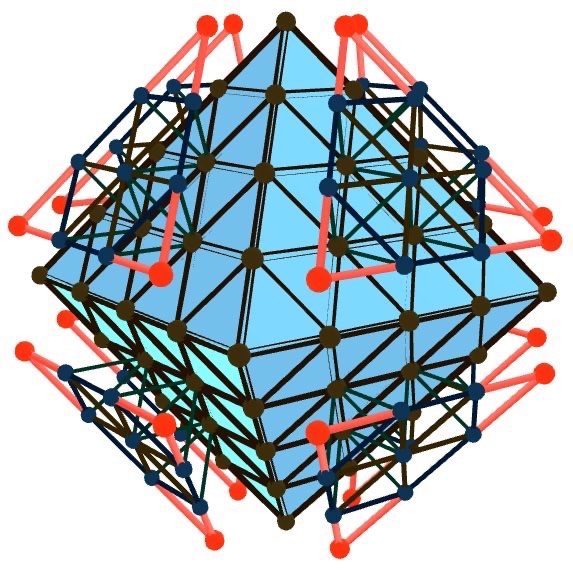}
  \caption*{Step 2}
\end{minipage}
\hfill
\begin{minipage}[b]{0.31\textwidth}
  \centering
  \includegraphics[width=\textwidth]{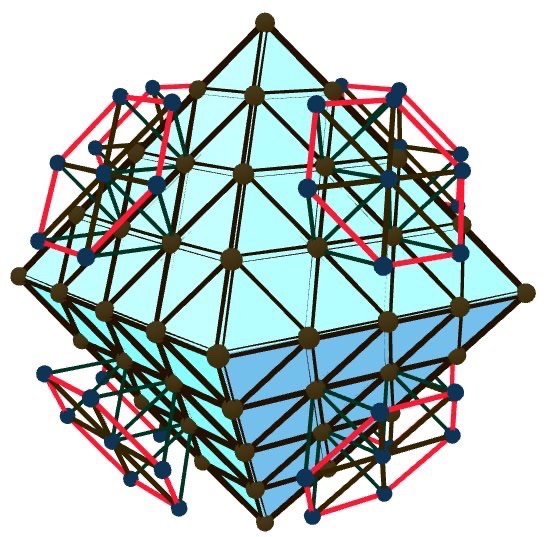}
  \caption*{Step 3}
\end{minipage}

\vspace{0.5cm}

\begin{minipage}[b]{0.31\textwidth}
  \centering
  \includegraphics[width=\textwidth]{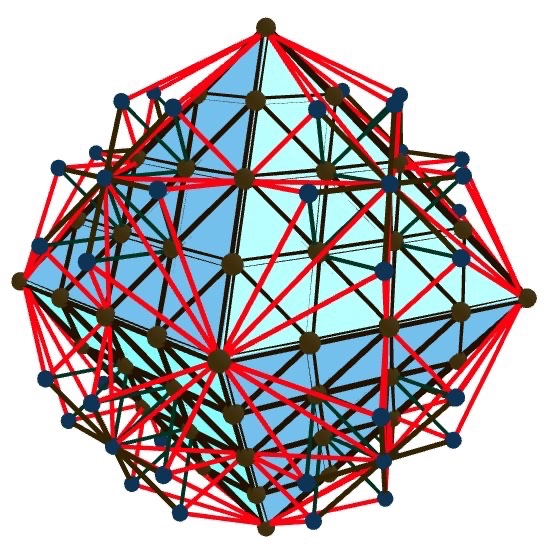}
  \caption*{Step 4}
\end{minipage}
\hfill
\begin{minipage}[b]{0.31\textwidth}
  \centering
  \includegraphics[width=\textwidth]{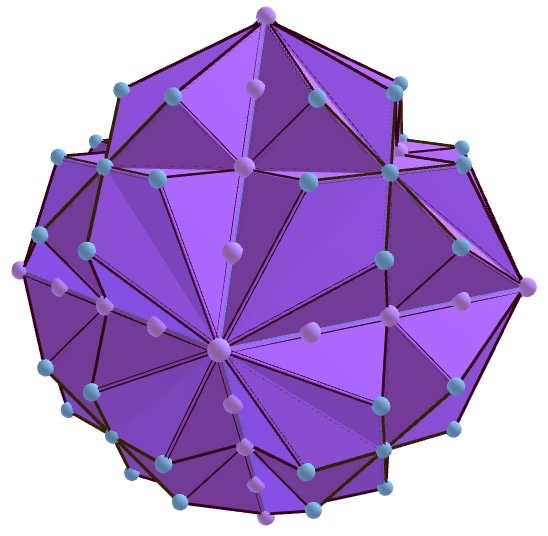}
  \caption*{Step 5}
\end{minipage}
\caption{The five construction steps of the M144p. See the main text for details.}
\label{fig:M144pSteps}
\end{figure}

\begin{enumerate}
\item \textbf{Start with a 6-frequency octahedron.}
  This is a subdivision of the regular octahedron in the FCC lattice
  with 6 steps per edge. Its interior consists of an arrangement of
  small tetrahedral and octahedral cells, characteristic of the FCC
  lattice. A 4-frequency octahedron sits concentrically inside the
  6-frequency octahedron. Both are visible in Step~1 of
  Figure~\ref{fig:M144pSteps}.

\item \textbf{Remove all nodes along the 12 edges of the 6-frequency
  octahedron.}
  What remains of the 6-frequency structure are nodes forming
  triangles over each face of the interior 4-frequency octahedron,
  as shown in Step~2.

\item \textbf{Remove the 3 corner vertices of each face triangle.}
  This leaves a network of nodes forming a hexagonal ring over each
  face of the 4-frequency octahedron, as shown in Step~3.

\item \textbf{Remove the hexagons' outer edges and make new connections.}
  The face-centre node of each 6-frequency face is connected to the
  6 exposed vertices of the 4-frequency octahedron lying on that face's
  boundary, and also to the mid-edge nodes of the 4-frequency octahedron.
  The face-centre nodes are already connected to the 6 remaining
  hexagon nodes, completing a star-like pattern on each face,
  as shown in Step~4.

\item \textbf{Fill in the triangular faces.}
  The exposed node network is filled in resulting in
  144 triangular faces of the M144p, as shown in Step~5.
\end{enumerate}

The result is the M144p: a closed triangulated surface with 74 vertices,
216 edges of various lengths,
and 144 triangular faces, with Euler characteristic $V - E + F = 2$.

\subsection{The inner icosahedron}
\label{sec:innerico}

When the 62 M120p vertices were matched to elements of the binary
icosahedral group $2I$ (Section~\ref{sec:lift}), 24 elements remained
unaccounted for---12 in the upper hemisphere at $w = +\varphi/2$ and
12 in the lower hemisphere at $w = -\varphi/2$. These are the
``missing 12.''

Under stereographic projection, the 12 upper-hemisphere elements
project to an icosahedron deep inside the M120p, at a stereographic
radius of $1/\sqrt{4\varphi+3}$---closer to the origin than any M120p
vertex. The directions of these 12 vertices are exactly the same as
the 12 C-vertex directions: cyclic permutations of
$(0, \pm\varphi, \pm\varphi^2) = \varphi \times (0, \pm 1, \pm\varphi)$.
Every inner icosahedron vertex is radially aligned with its
corresponding C-vertex (dot product $= 1.000$ for all 12 pairs), at a
radius ratio of exactly $1/\varphi$.

The inner icosahedron thus marks the deepest reach of the 600-cell
into the interior of the Mereon System. Its 12 vertices sit on the
focusing sphere (Section~\ref{sec:fs}), giving that sphere its first
geometric structure from the 600-cell correspondence. The full
derivation is given in Sections~\ref{sec:unaccounted} and~\ref{sec:inner}.

\subsection{The focusing sphere}
\label{sec:fs}

Between M144p and M120p lies a sphere with no polyhedral structure,
first described by Dennis \cite{Dennis1997}. An inner icosahedron,
whose 12 vertices are elements of $2I$ at latitude $|w| = \varphi/2$
on $S^3$ (Section~\ref{sec:inner}), has its vertices coincident with
this sphere. The function
of the sphere in the dynamics of the system is discussed in references~\cite{Dennis1997} and~\cite{GrayDennisKauffman2018}.

\subsection{The trefoil knot}
\label{sec:knot}

The Mereon trefoil is a $T(3,2)$ torus knot---the simplest non-trivial
knot, where $T(p,q)$ denotes a curve winding $p$ times around the
longitude and $q$ times through the meridian of the torus
\cite{Gray2002Knot, GrayHelixKnot, GrayKnotTypes, GrayOctaKnot, KauffmanKnotDynamics}.
The exponents $(2,3,5)$ of the Brieskorn singularity encode the
trefoil as the $(2,3)$ torus knot (Section~\ref{sec:eigenform}),
connecting the knot directly to the icosahedral symmetry of the Mereon System.

\begin{figure}[H]
\centering
\begin{minipage}[b]{0.44\textwidth}
  \centering
  \includegraphics[width=\textwidth]{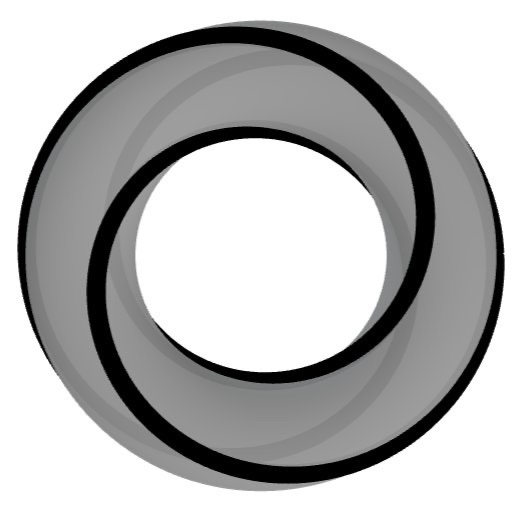}
  \caption*{\textbf{Mereon knot $T(3,2)$}}
\end{minipage}
\hfill
\begin{minipage}[b]{0.44\textwidth}
  \centering
  \includegraphics[width=\textwidth]{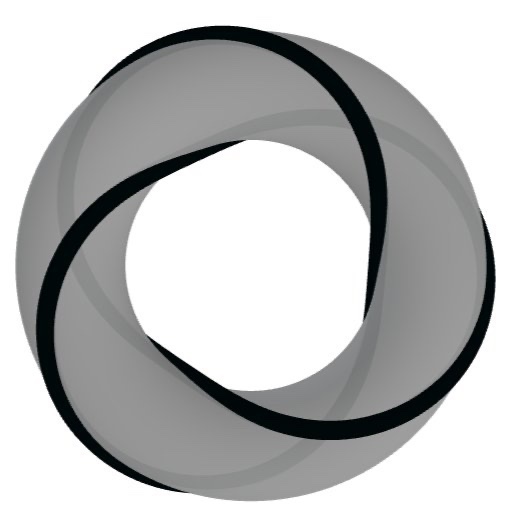}
  \caption*{\textbf{Standard trefoil $T(2,3)$}}
\end{minipage}
\caption{The two torus knot conformations of the trefoil, shown on their ring torus.
The Mereon knot $T(3,2)$ (left) winds three times around the longitude and twice
through the meridian, giving it 2-fold symmetry when viewed face-on.
The standard trefoil $T(2,3)$ (right) winds twice around the longitude and three
times through the meridian, giving it 3-fold symmetry face-on.
Both curves are topologically identical (same knot invariants) but wind
differently on the torus surface.}
\label{fig:torusknots}
\end{figure}

The trefoil appears in two geometric conformations, $T(3,2)$ and $T(2,3)$,
which wind differently on a torus surface but are topologically identical
(Figure~\ref{fig:torusknots}). The full treatment of these two conformations,
their relationship via Clifford torus projection in $S^3$, and their role as
the branch set of $M(2,3,5)$ is given in Section~\ref{sec:clifford}.

\subsection{The M120p boundary: non-crystallographic ($H_3$)}
\label{sec:boundary}

The M120p is a mixed convex-concave polyhedron with icosahedral
symmetry. Its symmetry group is the full icosahedral group $I_h$,
the Coxeter group $H_3$ (order 120), whose rotation subgroup
$I \cong A_5$ has order 60.

\begin{figure}[H]
\centering
\begin{minipage}[b]{0.31\textwidth}
  \centering
  \includegraphics[width=\textwidth]{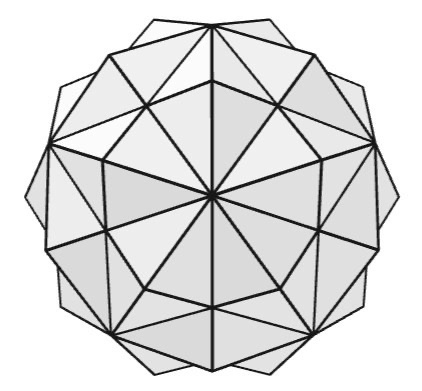}
\end{minipage}
\hfill
\begin{minipage}[b]{0.31\textwidth}
  \centering
  \includegraphics[width=\textwidth]{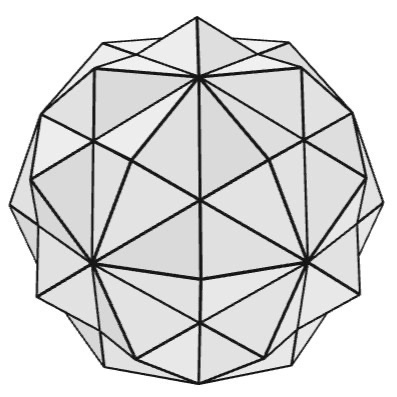}
\end{minipage}
\hfill
\begin{minipage}[b]{0.31\textwidth}
  \centering
  \includegraphics[width=\textwidth]{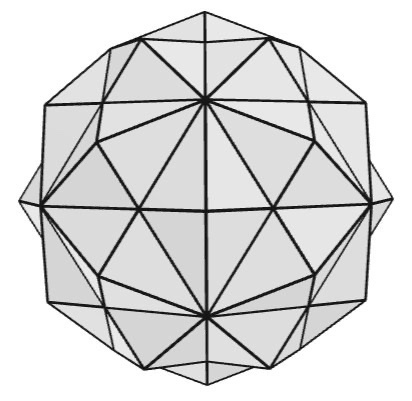}
\end{minipage}
\caption{The M120p 62-vertex boundary polyhedron (three views).}
\label{fig:M120p}
\end{figure}

\noindent The table below lists the properties of each vertex type.
The \textbf{Radius} column gives the Euclidean distance from the
origin to each vertex in Gray's coordinate system.
The \textbf{Stabiliser} column gives the subgroup of $I_h$ that fixes a
representative vertex---that is, the symmetries of the full icosahedral
group that leave that particular vertex unmoved. The notation $C_n$
denotes the cyclic group of order $n$: a $C_3$ stabiliser means a
vertex sits on a 3-fold rotation axis (unchanged by $120^\circ$ rotations),
$C_5$ means a 5-fold axis ($72^\circ$ rotations), and $C_2$ means a
2-fold axis ($180^\circ$ rotations). The \textbf{Fold} column repeats
this rotation order for clarity.

\begin{center}
\begin{tabular}{llllll}
\toprule
Type & Geometry & Count & Radius & Stabiliser & Fold \\
\midrule
A & Dodecahedron & 20 & $\varphi^2\sqrt{3} \approx 4.535$ & $C_3$ & 3 \\
C & Icosahedron & 12 & $\varphi^2\sqrt{1+\varphi^2} \approx 4.980$ & $C_5$ & 5 \\
B & Icosidodecahedron & 30 & $2\varphi^2 \approx 5.236$ & $C_2$ & 2 \\
\bottomrule
\end{tabular}
\end{center}

\noindent The three types do not sit on a single sphere. Each occupies
its own shell: A innermost, C intermediate, B outermost.

Each of the 120 faces contains one vertex of each type (the
\emph{trinity}).

The boundary (M120p) lives in $\mathbb{Q}(\varphi) = \mathbb{Q}(\sqrt{5})$---a
non-crystallographic number field, algebraically disjoint from the
core's integer coordinates.

\subsubsection*{Construction of the M120p}

The M120p is constructed in five steps, illustrated in
Figure~\ref{fig:M120pSteps}. The vertex coordinates of the
completed polyhedron are given in Appendix~\ref{app:matchtable} (the M120p column of the correspondence table).

\textbf{Step 1: The dodecahedron and the compound of five cubes.}
Begin with a regular dodecahedron, whose 20 vertices become the
Type~A vertices of the M120p. The dodecahedron contains within it
a compound of five cubes: five cubes share the dodecahedron's
20 vertices, each cube contributing 8 of the 20 vertices.

\textbf{Step 2: Adding the five dual octahedra.}
Each of the five cubes has a dual polyhedron, the octahedron.
Five octahedra are added to the model, scaled so that the
mid-edge points of each octahedron coincide with the mid-edge
points of its dual cube. The 30 vertices of these five octahedra
(6 vertices each) become the Type~B vertices of the M120p.
The compound of five octahedra in this step contributes 30 vertices;
connecting these 30 vertices to each other would form the
icosidodecahedron. For the M120p construction, these 30 vertices
are not connected to each other.

\textbf{Step 3: Adding the dual icosahedron.}
The dual of the dodecahedron is the icosahedron. An icosahedron
is added to the model, scaled so that its mid-edge points
coincide with the mid-edge points of the dodecahedron.
The 12 vertices of the icosahedron become the Type~C vertices
of the M120p.

\textbf{Step 4: The Rhombic Triacontahedron.}
The dodecahedron and the icosahedron, at this relative scaling,
together define the Rhombic Triacontahedron: the 20 dodecahedral
and 12 icosahedral vertices are precisely the vertices of
this well-known polyhedron.

\textbf{Step 5: Completing the M120p.}
The 30 Type~B vertices (from the five octahedra) are connected
to the 32 vertices of the Rhombic Triacontahedron (20 Type~A
plus 12 Type~C). This produces the 120 triangular faces of the
M120p, with each face containing exactly one vertex of each
type (the trinity A, B, C). The total vertex count is
$20 + 12 + 30 = 62$.

\begin{figure}[H]
\centering
\begin{minipage}[b]{0.31\textwidth}
  \centering
  \includegraphics[width=\textwidth]{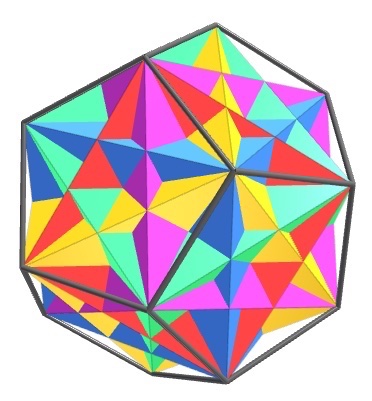}
  \caption*{Step 1}
\end{minipage}
\hfill
\begin{minipage}[b]{0.31\textwidth}
  \centering
  \includegraphics[width=\textwidth]{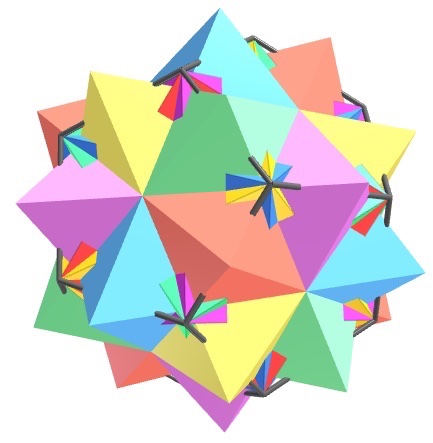}
  \caption*{Step 2}
\end{minipage}
\hfill
\begin{minipage}[b]{0.31\textwidth}
  \centering
  \includegraphics[width=\textwidth]{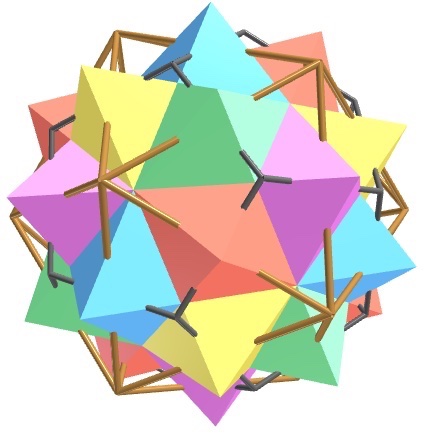}
  \caption*{Step 3}
\end{minipage}

\vspace{0.5cm}

\begin{minipage}[b]{0.31\textwidth}
  \centering
  \includegraphics[width=\textwidth]{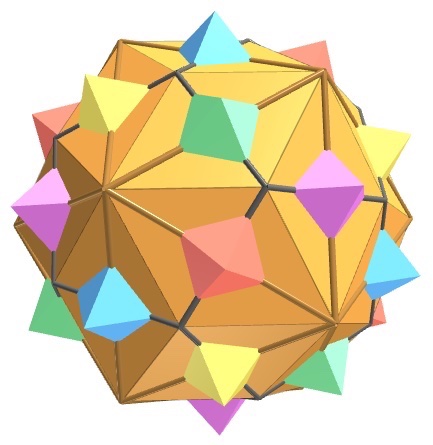}
  \caption*{Step 4}
\end{minipage}
\hfill
\begin{minipage}[b]{0.31\textwidth}
  \centering
  \includegraphics[width=\textwidth]{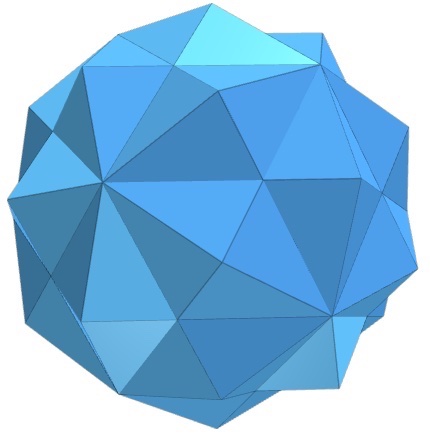}
  \caption*{Step 5}
\end{minipage}
\caption{The five construction steps of the M120p. See the main text for details.}
\label{fig:M120pSteps}
\end{figure}

\subsection{Crystallographic meets non-crystallographic}
\label{sec:meeting}

The Mereon System is the meeting of two worlds:

\begin{center}
\begin{tabular}{lll}
\toprule
& M144p core & M120p boundary \\
\midrule
Symmetry & $O_h$ (crystallographic) & $I_h$ (non-crystallographic, $H_3$) \\
Number field & $\mathbb{Q}$ (integers) & $\mathbb{Q}(\varphi) = \mathbb{Q}(\sqrt{5})$ \\
Lattice & FCC & None (cannot tile space) \\
McKay & $E_7$ (via $2O$) & $E_8$ (via $2I$) \\
\bottomrule
\end{tabular}
\end{center}

\noindent The Mereon System unites these two algebraically disjoint symmetry
families by \emph{spatial nesting}: the $O_h$-symmetric core sits
inside the $H_3$-symmetric boundary, with the focusing sphere marking
the boundary between the crystallographic and non-crystallographic worlds.

\section{$H_3$, $H_4$, and the 600-Cell}
\label{sec:H3H4}

The Coxeter group $H_4$ is the symmetry group of the 600-cell (and its
dual, the 120-cell) in $\mathbb{R}^4$. It has order 14400. Its
rotational subgroup has order 7200. $H_4$ is the \emph{unique}
non-crystallographic Coxeter group in four dimensions \cite{Humphreys1990}.

The inclusion $H_3 \subset H_4$ is realised by the embedding
$\mathbb{R}^3 \hookrightarrow \mathbb{R}^4$: any 3-dimensional
hyperplane section of the 600-cell inherits (at most) $H_3$ symmetry,
and a section through the equator ($w = 0$) achieves it exactly.

The binary icosahedral group $2I \subset \mathrm{SU}(2) \cong S^3$ has
order 120. Its elements, as unit quaternions, are the 120 vertices of
the 600-cell. We recall them:

\begin{enumerate}[label=(\alph*)]
\item 8 elements: $\pm 1, \pm i, \pm j, \pm k$.
\item 16 elements: $\tfrac{1}{2}(\pm 1 \pm i \pm j \pm k)$.
\item 96 elements:
  $\tfrac{1}{2}(0, \pm 1, \pm\varphi, \pm\varphi^{-1})$ and even
  permutations.
\end{enumerate}

\noindent Families (a) and (b) form $2T$ (24 elements, binary
tetrahedral). The full 120 form $2I$ (binary icosahedral).

The 24 Hurwitz units in $2T$ have coordinates in $\mathbb{Q}$---no
golden ratio required. The 96 golden quaternions in $2I \setminus 2T$
have coordinates in $\mathbb{Q}(\varphi)$. Without $\varphi$, these
quaternions cannot be written down. They represent the additional
symmetries that distinguish icosahedral from tetrahedral symmetry.

\subsection{The Icosians and the $E_8$ Lattice}
\label{sec:icosians}

The integer linear combinations of the 120 elements of $2I$ form a
subring of the quaternions called the \textbf{icosians}
\cite{ConwaySloane1999, Baez2017}. Every icosian has components in
$\mathbb{Q}(\sqrt{5})$, so icosians can be represented as 8-tuples of
rationals forming a lattice in $\mathbb{R}^8$. Equipping this lattice
with the modified norm $\|q\|^2 = x + y$ (where $|q|^2 = x +
\sqrt{5}\,y$) yields the $E_8$ lattice \cite{ConwaySloane1999,Baez2017}.

Conway and Sloane \cite{ConwaySloane1999} and Baez \cite{Baez2017} established
that the icosian ring, equipped with the modified norm $\|q\|^2 = x + y$
(where $|q|^2 = x + \sqrt{5}\,y$), is isometrically isomorphic to the
$E_8$ lattice in $\mathbb{R}^8$.

\subsection{The Du Val Resolution and the $E_8$ Dynkin Diagram}
\label{sec:duval}

A second route from the icosahedron to $E_8$ goes through algebraic
geometry \cite{duVal1934, Baez2017}. Klein \cite{Klein1888} showed
that every $2I$-invariant polynomial on $\mathbb{C}^2$ is a
polynomial in three invariants $V$, $E$, $F$ satisfying
$V^5 + E^2 + F^3 = 0$, so $\mathbb{C}^2/2I \cong V(2,3,5)$ — the
Brieskorn variety central to this paper. The minimal resolution of
this singularity replaces the origin with 8 copies of $\mathbb{CP}^1$
intersecting in the $E_8$ Dynkin diagram pattern
\cite{duVal1934, Baez2017, KirbyScharlemann1982}:

\begin{center}
\begin{tikzpicture}[scale=0.8]
  \foreach \x in {1,...,7}
    \fill (\x,0) circle (3pt);
  \fill (4,1) circle (3pt);
  \foreach \x in {1,...,6}
    \draw (\x,0) -- (\x+1,0);
  \draw (4,0) -- (4,1);
\end{tikzpicture}
\end{center}

As Baez \cite{Baez2017} notes, the icosian construction and the du~Val
resolution both derive $E_8$ from $2I$ by entirely different routes;
their connection remains open. The eigenform loop
(Section~\ref{sec:eigenform}) provides structural context suggesting
such a connection exists. The Poincar\'{e} homology sphere is
$M(2,3,5) = \mathrm{SU}(2)/2I$, the subspace of unit-size icosahedra
in $\mathbb{C}^2/2I$ \cite{Baez2017}.

\section{The Drop-$w$ Lift: From 3D Coordinates to Unit Quaternions}
\label{sec:lift}

This section describes how the 600-cell correspondence was discovered.

\subsection{The scaling}

A unit quaternion $q = (w, x, y, z)$ satisfies $w^2 + x^2 + y^2 + z^2 = 1$.
Dropping $w$ gives a 3D point $(x, y, z)$ with $r^2 = x^2 + y^2 + z^2 = 1 - w^2 \leq 1$.
But Gray's M120p coordinates have $r_{\max} = 2\varphi^2 \approx 5.236$
(the B-vertices at the outermost shell).

The B-vertices at $r = 2\varphi^2$ are the outermost. If these correspond
to the equatorial elements of $2I$ (those with $w = 0$, hence $r = 1$ on $S^3$),
the scale factor is $2\varphi^2$. We therefore normalise:
\[
  r' = (x', y', z') = \frac{(x, y, z)}{2\varphi^2}.
\]

After scaling, the table below shows the resulting scaled radial
distance $r'$ for each vertex type. Here $r'$ is the Euclidean
distance from the origin after dividing by $2\varphi^2$, so that
the outermost B-vertices land exactly at $r' = 1$ on the unit sphere
$S^3$. The values of $r'$ are what the $x$, $y$, $z$ components of a
unit quaternion can reach when $w \geq 0$:
\begin{center}
\begin{tabular}{lll}
\toprule
Type & Scaled radius $r'$ & Value \\
\midrule
B & $1$ & $= 1$ (equator) \\
A & $\sqrt{3}/2$ & $\approx 0.866$ \\
C & $\sin(72^\circ)$ & $\approx 0.951$ \\
\bottomrule
\end{tabular}
\end{center}

\subsection{Recovering $w$}

Since $r'^2 + w^2 = 1$ on the unit $S^3$:
\begin{equation}\label{eq:wrecovery}
  w = +\sqrt{1 - r'^2}.
\end{equation}

The positive sign selects the upper (northern) hemisphere. The recovered
$w$-values are shown in the table below. Each quaternion $q \in 2I$
represents a rotation of $\mathbb{R}^3$ by an angle $\theta$, where
$w = \cos(\theta/2)$; the \textbf{Rotation angle} column gives this
angle $\theta$, which is the angle through which the corresponding
element of the icosahedral rotation group $I$ rotates space:
\begin{center}
\begin{tabular}{llll}
\toprule
Type & $w$ & Exact value & Rotation angle \\
\midrule
A & $1/2$ & $= \cos(60^\circ)$ & $120^\circ$ \\
C & $1/(2\varphi)$ & $= \cos(72^\circ)$ & $144^\circ$ \\
B & $0$ & $= \cos(90^\circ)$ & $180^\circ$ \\
\bottomrule
\end{tabular}
\end{center}

\subsection{The 62 of 62 vertex match}

The M120p vertices correspond exactly to 62 of the 120 elements of $2I$
via the following four-step procedure. First, each M120p vertex $(x,y,z)$
is divided by $2\varphi^2$ to obtain a scaled point $(x',y',z')$ with
$r' = |(x',y',z')| \leq 1$. Second, the missing fourth coordinate is
recovered by setting $w = +\sqrt{1 - r'^2}$, placing the point on the
upper hemisphere of $S^3$. Third, the result $(w, x', y', z')$
automatically satisfies $w^2 + x'^2 + y'^2 + z'^2 = 1$ --- it is a
unit quaternion requiring no further scaling. Fourth, and crucially,
each of the 62 resulting unit quaternions coincides exactly with one of
the 120 elements of $2I$. The full explicit correspondence, vertex by
vertex, is tabulated in Appendix~\ref{app:matchtable}, where every entry
in the ``Lifted'' and ``Matched $2I$ element'' columns is identical.

The 62 matched elements are precisely those in the upper hemisphere
($w \geq 0$) of $2I$ with $w \in \{0,\; \tfrac{1}{2},\; \tfrac{1}{2\varphi}\}$,
that is, the 20 Type~A vertices ($w = \tfrac{1}{2}$), the 12 Type~C
vertices ($w = \tfrac{1}{2\varphi}$), and the 30 Type~B vertices ($w = 0$).
The identity element $q = +1$ and the 12 upper-hemisphere elements at
$w = \varphi/2$ are not M120p vertices but belong to the Mereon System:
the 12 at $w = \varphi/2$ form the inner icosahedron at the focusing
sphere (Section~\ref{sec:unaccounted}).

The remaining 58 elements of $2I$ --- those in the lower hemisphere ($w < 0$),
the identity and anti-identity, and the 12 upper-hemisphere elements
at $w = \varphi/2$ --- do not appear as M120p vertices. The
lower-hemisphere elements form the reciprocal Mereon System
(Section~\ref{sec:shells}).

\subsection{The remaining 24 elements and the inner icosahedron}
\label{sec:unaccounted}

The 120 elements of $2I$ distribute across five $|w|$-values:

\begin{center}
\begin{tabular}{lllll}
\toprule
$|w|$ & Upper ($w > 0$) & Equator ($w = 0$) & Lower ($w < 0$) & Total \\
\midrule
$1$ (identity) & 1 & --- & 1 & 2 \\
$\varphi/2$ ($36^\circ$ lat.) & 12 & --- & 12 & 24 \\
$1/2$ ($60^\circ$ lat.) & 20 & --- & 20 & 40 \\
$1/(2\varphi)$ ($72^\circ$ lat.) & 12 & --- & 12 & 24 \\
$0$ (equator) & --- & 30 & --- & 30 \\
\bottomrule
\end{tabular}
\end{center}

\noindent The total element count is: $2 + 24 + 40 + 24 + 30 = 120$.

The M120p uses 62 elements, all from the upper hemisphere ($w \geq 0$):
the 20 with $w = +1/2$ (A-type), the 12 with $w = +1/(2\varphi)$
(C-type), and the 30 with $w = 0$ (B-type, the equatorial elements).
The corresponding negative-$w$ elements---the 20 with $w = -1/2$ and
the 12 with $w = -1/(2\varphi)$ in the lower hemisphere---are the
``hidden'' portion of the double cover, with no counterpart among
M120p vertices. The 24 elements at
$|w| = \varphi/2$ (latitude $36^\circ$) sit at the highest latitude
of any non-identity element in $2I$:
$w = \varphi/2$, exactly $36^\circ$ from the north pole of $S^3$.
In the stereographic projection, they map to Shell~1 at radius
$r = 1/\sqrt{4\varphi+3} \approx 0.3249$---the
innermost non-trivial shell, closer to the origin than any M120p
vertex type. These 12 vertices of Shell~1 form an inner icosahedron.
We assign these 12 vertices to coincide with the focusing sphere
in the Mereon System (Section~\ref{sec:inner}).

\section{The $H_3 \subset H_4$ Projection}
\label{sec:main}

\subsection{Angular alignment of the 600-cell projection with the M120p}

The stereographic projection $\pi$ defined in equation~\eqref{eq:stereo}
maps every element $q \in 2I$ with $q \neq \pm 1$ to a point in
$\mathbb{R}^3$. The key geometric fact is that every such projected
point lies \emph{exactly} along one of the 62 M120p vertex directions.
This is not a numerical approximation: it is an exact consequence of
the structure of $2I$.

To see why, write any element of $2I$ in the form
$q = \cos(\theta/2) + \sin(\theta/2)\hat{n}$, where $\theta$ is a
rotation angle and $\hat{n}$ is the unit rotation axis. The imaginary
part of $q$ is $\mathrm{Im}(q) = \sin(\theta/2)\hat{n}$, which points
along $\hat{n}$. The stereographic projection of $q$ is
\[
  \pi(q) = \frac{(x,y,z)}{1+w},
\]
whose direction in $\mathbb{R}^3$ is determined solely by
$\mathrm{Im}(q)$ and hence by the rotation axis $\hat{n}$.
The rotation axes of the icosahedral group $I$ are exactly the 62
M120p vertex directions: the 20 dodecahedral axes (Type~A, 3-fold),
the 12 icosahedral axes (Type~C, 5-fold), and the 30
icosidodecahedral axes (Type~B, 2-fold). Since every non-trivial
element of $2I$ rotates $\mathbb{R}^3$ about one of these axes, every
projected 600-cell vertex points exactly along an M120p vertex
direction. This is the $H_3 \subset H_4$ correspondence made explicit:
the non-crystallographic symmetry of the M120p in $\mathbb{R}^3$ is
the restriction of the full $H_4$ symmetry of the 600-cell in
$\mathbb{R}^4$, and stereographic projection is the map that makes
this restriction precise.

\subsection{The eight-shell structure}
\label{sec:shells}

The distance from the origin after stereographic projection depends
only on the scalar part $w$ of the quaternion:
\[
  r = |\pi(q)| = \sqrt{\frac{1-w}{1+w}}.
\]
Since the scalar parts of $2I$ take only the values
$\{0, \pm\tfrac{1}{2}, \pm\tfrac{\varphi}{2}, \pm\tfrac{1}{2\varphi}, \pm 1\}$,
the 120 projected elements fall into exactly 8 shells (plus the origin
$r=0$ and the point at infinity).

Each $|w|$-value determines both the shell radius and the type of
rotation axis, since a quaternion with scalar part $w = \cos(\theta/2)$
represents a rotation by angle $\theta$. The correspondence is:

\begin{center}
\begin{tabular}{llll}
\toprule
$|w|$ & Rotation angle $\theta$ & Axes & Type \\
\midrule
$\varphi/2$ & $72^\circ$ & Icosahedron vertices & C \\
$1/2$ & $120^\circ$ & Dodecahedron vertices & A \\
$0$ & $180^\circ$ & Edge midpoints & B \\
$1/(2\varphi)$ & $144^\circ$ & Icosahedron vertices & C \\
\bottomrule
\end{tabular}
\end{center}

\noindent\textbf{Note on the Type column in the following table.} ``Type'' denotes the
\emph{axis direction} of the rotation: A = dodecahedral (3-fold axes),
B = edge-midpoint (2-fold axes), C = icosahedral (5-fold axes). Both
shells~1 and~3 (and their reciprocals 7 and~5) point along icosahedral
axes, hence both are labeled C. However, only shells~3 and~5 correspond to the M120p axes defined by
C-vertices ($w = \pm\tfrac{1}{2\varphi}$, 144° rotations).
Shells~1 and~7 ($w = \pm\tfrac{\varphi}{2}$, 72° rotations) point
along the same 12 icosahedral axes and correspond to the 24 elements
that form the inner icosahedron (Section~\ref{sec:unaccounted})
at the focusing sphere.

The complete shell table, with exact algebraic expressions for each radius,
is as follows. The expressions are derived from the golden-ratio identities
$2\varphi - 1 = \sqrt{5}$, $2\varphi + 1 = \varphi^3$, and
$4\varphi + 3 = \varphi^2(\varphi^2+1)$; full derivations are given
in Section~\ref{sec:shellderivations} below.

\begin{center}
\begin{tabular}{lllll}
\toprule
Shell & $r$ & Expression & Type & Count \\
\midrule
0 & $0$ & $q = +1$ & --- & 1 \\
1 & $0.3249$ & $1/\sqrt{4\varphi+3}$ & C & 12 \\
2 & $0.5774$ & $1/\sqrt{3}$ & A & 20 \\
3 & $0.7265$ & $5^{1/4}/\varphi^{3/2}$ & C & 12 \\
4 & $1.0000$ & $1$ & B & 30 \\
5 & $1.3764$ & $\varphi^{3/2}/5^{1/4}$ & C & 12 \\
6 & $1.7321$ & $\sqrt{3}$ & A & 20 \\
7 & $3.0777$ & $\sqrt{4\varphi+3}$ & C & 12 \\
$\infty$ & $\infty$ & $q = -1$ & --- & 1 \\
\bottomrule
\end{tabular}
\end{center}

\noindent $1 + 12 + 20 + 12 + 30 + 12 + 20 + 12 + 1 = 120$. Shell
counts per type: $|A| = 2$, $|B| = 1$, $|C| = 4$.

The M120p vertices occupy the upper hemisphere ($w \geq 0$): shells~2, 3,
and~4, giving the 20 Type~A, 12 Type~C, and 30 Type~B vertices
respectively. Shell~1 contains the 12 inner icosahedron vertices at the
focusing sphere, and Shell~0 is the north pole ($q = +1$).
Together, shells~0 through~4 constitute the Mereon System as
projected into $\mathbb{R}^3$.

Shells~4 through~$\infty$ form the reciprocal Mereon System: the
lower-hemisphere counterpart. For every vertex at stereographic
distance $r$ in the upper hemisphere, there is a corresponding vertex
at distance $1/r$ in the lower hemisphere, with $r \leq 1$ for the
upper-hemisphere shells (after the inner icosahedron). Shell~4 (the
equatorial B-vertices at $r = 1$) is shared between both halves.
The full 120-element structure of $2I$ is thus partitioned into the
Mereon System and its reciprocal, with no elements left over.

The shells come in reciprocal pairs: $|\pi(q)| \cdot |\pi(-q)| = 1$,
since the product $\sqrt{\frac{1-w}{1+w} \cdot \frac{1+w}{1-w}} = 1$.
This reciprocal pairing is the visible signature of the
$\{\pm 1\}$ kernel of the double cover $\mathrm{SU}(2) \to \mathrm{SO}(3)$:
$q$ and $-q$ represent the same rotation but land at reciprocal
distances under stereographic projection.
Type~A has 2 shells (at $r = 1/\sqrt{3}$ and $\sqrt{3}$, a reciprocal pair);
Type~B has 1 self-reciprocal shell at $r = 1$;
Type~C has 4 shells forming two reciprocal pairs.

The 30 Type~B projections at $r = 1$ are the pure imaginary quaternions
in $2I$: the half-turns of $I$, sitting exactly at the equator $w = 0$ of $S^3$.

The type-pure shell structure shows that the $H_3$ decomposition of
$\mathbb{R}^3$ into three orbit types (A, B, C vertex directions)
lifts to an $H_4$-compatible decomposition of the 600-cell into
conjugacy classes. The orbit structure is completely preserved by
the stereographic projection.

\subsection{Two correspondences: vertices and faces}

The 62 of 62 vertex match of Section~\ref{sec:lift} is one correspondence
between the M120p and $2I$. There is a second, entirely separate
correspondence involving the \emph{faces}.

The 120 face centroids of the M120p all lie on a single sphere of
radius $\approx 4.6950$, because the icosahedral group $I_h$ acts
transitively on the faces: it can map any face to any other face by a
symmetry transformation. Since there are 120 faces and $|2I| = 120$,
there is a natural one-to-one correspondence between the elements of
$2I$ and the faces of the M120p.

More precisely: a group $G$ acts \emph{transitively} on a set $X$ when
for any two elements $x, y \in X$ there exists a group element $g \in G$
with $g \cdot x = y$---that is, every element of $X$ is reachable from
every other by some symmetry. Here, $I_h$ acts transitively on the
120 triangular faces of the M120p, meaning every face is equivalent to
every other face under the icosahedral symmetry. Since $|I_h| = 120$
equals the number of faces, each element of $I_h$ (and hence each element
of its double cover $2I$) corresponds to exactly one face and vice versa.

This face correspondence is fundamentally different from the vertex
correspondence. In the vertex match, 62 elements of $2I$ (those in the
upper hemisphere, $w \geq 0$) are matched to the 62 vertices via the
directions of stereographic projection. In the face correspondence, all
120 elements of $2I$ are matched to the 120 faces via the transitive
group action, with no projection involved. The two correspondences
reflect two distinct mathematical relationships between $2I$ and the M120p.

\subsection{Derivation of shell radii}
\label{sec:shellderivations}

The 8 shells arise because the stereographic projection radius
$r = \sqrt{(1-w)/(1+w)}$ maps each $|w|$-value to two shells:
one for $w > 0$ (upper hemisphere) and one for $-w < 0$ (lower
hemisphere), forming the reciprocal pairs $r(w) \cdot r(-w) = 1$.
The M120p vertices occupy shells 2, 3, and 4 (types A, C, B).
Shell~1 contains the inner icosahedron at the focusing sphere.
The remaining shells (5, 6, 7) are the reciprocal counterparts
in the lower hemisphere (Section~\ref{sec:unaccounted}).

\medskip
\noindent\textbf{Shell 0} ($w = +1$, identity):
$r = \sqrt{0/2} = 0$.

\medskip
\noindent\textbf{Shell 1} ($w = +\varphi/2$, Type C):
\[
  r = \sqrt{\frac{2-\varphi}{2+\varphi}} = \frac{1}{\sqrt{4\varphi+3}},
\]
since $(2-\varphi)(4\varphi+3) = 8\varphi+6-4\varphi^2-3\varphi = 5\varphi+6-4(\varphi+1)
= \varphi+2 = 2+\varphi$. Therefore
$(2-\varphi)/(2+\varphi) = 1/(4\varphi+3)$, giving
$r = 1/\sqrt{4\varphi+3} \approx 0.3249$.

\medskip
\noindent\textbf{Shell 2} ($w = +1/2$, Type A):
\[
  r = \sqrt{\frac{1/2}{3/2}} = \frac{1}{\sqrt{3}}.
\]

\medskip
\noindent\textbf{Shell 3} ($w = +1/(2\varphi)$, Type C):
\[
  r = \sqrt{\frac{2\varphi-1}{2\varphi+1}} = \sqrt{\frac{\sqrt{5}}{\varphi^3}}
    = \frac{5^{1/4}}{\varphi^{3/2}},
\]
using $2\varphi - 1 = \sqrt{5}$ and $2\varphi + 1 = \varphi^3$.

\medskip
\noindent\textbf{Shell 4} ($w = 0$, Type B):
$r = \sqrt{1/1} = 1$.

\medskip
\noindent\textbf{Shells 5, 6, 7} are the reciprocals of Shells 3, 2, 1
respectively, by the identity $r(w)\cdot r(-w) = 1$:
\[
  r_5 = \frac{\varphi^{3/2}}{5^{1/4}}, \qquad
  r_6 = \sqrt{3}, \qquad
  r_7 = \sqrt{4\varphi+3}.
\]

\subsection{The latitude partition}
\label{sec:latitude}

The recovered $w$-values reveal the structural meaning of the vertex
classification. Think of $S^3$ as a globe: the $w$-coordinate is
latitude. The vertex type is then determined entirely by latitude in the
fourth dimension. The classification that was discovered geometrically in
three dimensions turns out to be a latitude classification in four-dimensional
space.

The three M120p vertex types are:
\begin{itemize}
\item A ($w = 1/2$, latitude $60^\circ$): Input.
\item B ($w = 0$, equator, latitude $0^\circ$): Thruput.
\item C ($w = 1/(2\varphi)$, latitude $72^\circ$): Output.
\end{itemize}

\noindent Including the inner icosahedron at the focusing sphere
($|w| = \varphi/2$, latitude $36^\circ$), the non-equatorial latitudes
form a $\varphi$-ladder: each step outward from the inner icosahedron
multiplies $|w|$ by $1/\varphi$, so
$\varphi/2 \times 1/\varphi = 1/2$ (A-vertices) and
$1/2 \times 1/\varphi = 1/(2\varphi)$ (C-vertices). The shells
descend from the inner icosahedron to the equator by successive
divisions by the golden ratio.

\section{The Field Extension Interpretation}
\label{sec:failed}

It is well known that
$2O \not\subset 2I$: since $|2O| = 48$ does not divide $|2I| = 120$
(as $120/48 = 2.5 \notin \mathbb{Z}$), Lagrange's theorem immediately
excludes $2O$ as a subgroup of $2I$ \cite{Humphreys1990}.

The
cascade $3 \to 4 \to 5$ (the highest rotation orders of the tetrahedron,
octahedron, and icosahedron respectively) corresponds to successive field extensions:
\[
  \mathbb{Q} \;\longrightarrow\; \mathbb{Q}(\sqrt{2})
  \;\longrightarrow\; \mathbb{Q}(\sqrt{5}).
\]
This chain refers specifically to the quaternion components of the
binary groups as unit quaternions in $S^3$, not to the coordinate
fields of the corresponding polyhedra. The M144p can be scaled to have
integer coordinates in $\mathbb{Q}$ regardless of this chain; its
connection to $E_7$ is through the McKay correspondence applied
to its $O_h$ symmetry, not through its vertex coordinates.

The $2T$ elements (the 24 elements of $2T \subset 2I$) can be written with
rational quaternion coordinates alone. The $2O$ elements, as unit
quaternions, require $\sqrt{2}$ (e.g.\ $\tfrac{1}{\sqrt{2}}(\pm 1 \pm i)$).
The $2I$ elements require $\sqrt{5}$ (equivalently, $\varphi$). Each step
requires number-field resources the previous step cannot access.
In particular, every C-vertex (Output) requires $\varphi$ to be written
down as a quaternion coordinate. Without the golden ratio, the C-vertices cannot be realised. The cascade is
an expansion of what numbers are available.

The non-inclusion $2O \not\subset 2I$ means specifically that
Thruput~(B) does not sit inside Output~(C). The other two nestings hold:
\begin{itemize}
\item \textbf{$2T \subset 2O$:} every rotation of the tetrahedron is
  also a rotation of the octahedron, since $|T| = 12$ divides $|O| = 24$.
\item \textbf{$2T \subset 2I$:} every rotation of the tetrahedron is
  also a rotation of the icosahedron, since $|T| = 12$ divides $|I| = 60$.
\item \textbf{$2O \not\subset 2I$:} octahedral
  symmetries (4-fold) do \emph{not} extend to icosahedral symmetries
  (5-fold). The icosahedron has 2-fold, 3-fold, and 5-fold rotation
  axes, but no 4-fold axes. A $90^\circ$ rotation about a 4-fold axis
  of the octahedron simply does not exist as a symmetry of the
  icosahedron. The 4-fold world and the 5-fold world share the 3-fold
  core and nothing more.
\end{itemize}

Under the McKay correspondence \cite{McKay1980,Humphreys1990,Baez2017},
$2T \to E_6$, $2O \to E_7$, $2I \to E_8$. The genuine subgroup
inclusions ($2T \subset 2O$ and $2T \subset 2I$) are reflected
faithfully: $E_6$ appears as a sub-diagram in both $E_7$ and $E_8$.
But $2O \not\subset 2I$ means $E_7$ and $E_8$ are
\emph{siblings}---both containing $E_6$, but neither containing the
other. The structure is a Y, with $E_6$ at the junction and $E_7$, $E_8$
as the two branches (shown here laid out horizontally):
\[
  E_7 \;\longleftarrow\; E_6 \;\longrightarrow\; E_8
\]
with $E_6$ (the tetrahedral core, the shared $2T$) at the
junction. This Y-structure replaces the linear chain
$2T \subset 2O \subset 2I$.

The polyhedra coexist in $\mathbb{R}^3$ because geometric space
provides $\mathbb{Q}(\sqrt{2}, \sqrt{5})$ simultaneously. But $2O$
and $2I$, as algebraic objects, live in different number fields and
cannot be subgroups of one another. The focusing sphere bridges what
algebra alone cannot join.

\subsection{The 24-cell inside the 600-cell and the $E_6 \subset E_8$ embedding}

McKay's theorem \cite{McKay1980} states that the McKay graph of a
finite subgroup of $\mathrm{SU}(2)$ is the \emph{affine} Dynkin
diagram (denoted with a hat): $2T \to \hat{E}_6$, $2O \to
\hat{E}_7$, $2I \to \hat{E}_8$. The affine diagram $\hat{E}_n$
has one extra node compared to the ordinary diagram $E_n$.
Throughout the rest of this paper, $E_6$, $E_7$, $E_8$ refer to
the ordinary (non-affine) Lie algebras and root systems, with
dimensions 78, 133, and 248 respectively.

The 24 elements of $2T \subset 2I$, when stereographically projected,
fall into exactly three shells, forming the projected image of the 24-cell:

\begin{center}
\begin{tabular}{llll}
\toprule
$r$ & Expression & Count & Directions \\
\midrule
$1/\sqrt{3}$ & $0.577$ & 8 & Cube vertices (A-type) \\
$1$ & $1.000$ & 6 & Coordinate axes ($\pm i, \pm j, \pm k$) \\
$\sqrt{3}$ & $1.732$ & 8 & Cube vertices (A-type) \\
\bottomrule
\end{tabular}
\end{center}

\noindent Plus origin and infinity, giving 24 total. The 16 elements
with $|w| = 1/2$ project to A-type shells at $r = 1/\sqrt{3}$
(those with $w = +1/2$) and $r = \sqrt{3}$ (those with $w = -1/2$),
forming a reciprocal pair. The 6 elements $\pm i, \pm j, \pm k$ have
$w = 0$ and project to $r = 1$ along coordinate axes --- the 2-fold
axes shared by both octahedral and icosahedral symmetry. The 96 edges
of the 24-cell (at geodesic angle $\pi/3$ on $S^3$) are confirmed, and
all directions coincide exactly with tetrahedral/octahedral vertex
directions by the same angular alignment argument of Section~\ref{sec:main}.

The 600-cell contains exactly five inscribed 24-cells, each a distinct
conjugate copy of $2T$ within $2I$, with $5 \times 24 = 120$ vertices
accounting for the full 600-cell. This is a classical fact of $H_4$
geometry. It gives the $E_6 \subset E_8$ embedding: the 24-cell
($2T$, McKay $\to \hat{E}_6$) sits inside the 600-cell ($2I$, McKay
$\to \hat{E}_8$).

\subsection{The $E_7$ bridge: where crystallographic meets non-crystallographic}
\label{sec:E7}

The octahedral group $O$ has order 24 and the icosahedral group $I$
has order 60. Since 24 does not divide 60, $O$ cannot be a subgroup
of $I$, and likewise $2O \not\subset 2I$. McKay's theorem maps $2O$
to $E_7$ independently of the $2T \subset 2I$ embedding. The Mereon
System therefore realises the three exceptional algebras by two routes:

\begin{enumerate}
\item \textbf{Embedding:} $2T \subset 2I$ gives $E_6 \subset E_8$.
  This is the 24-cell inside the 600-cell. The geometry is $H_3 \subset H_4$.
\item \textbf{Independent correspondence:} $2O$ (order 48) gives $E_7$.
  The M144p core with its $O_h$ symmetry of order 48
  and its scalable-to-integer FCC coordinates is the $E_7$
  structure. It lives \emph{inside} the
  $E_8$ boundary physically but not algebraically.
\end{enumerate}

The $E_6 \subset E_8$ route is non-crystallographic: $H_3 \subset H_4$,
golden ratio throughout, no lattice. The $E_7$ route is crystallographic:
$O_h$, integer coordinates, FCC lattice. The Mereon System contains
\emph{both} routes in a single spatial structure.

The M144p, as an FCC lattice structure with $O_h$ symmetry, shares the
McKay correspondence $2O \to \hat{E}_7$ with any $O_h$-symmetric
crystallographic structure. Its significance here is not that it uniquely
realises $E_7$, but that it physically sits inside the M120p, which
realises $E_8$---placing an $E_7$-type structure spatially inside an
$E_8$-type structure within the Mereon System.

The focusing sphere marks a geometric boundary between two algebraically
distinct regions: the interior, governed by the crystallographic $O_h$
symmetry of the M144p and corresponding to $E_7$, and the exterior,
governed by the non-crystallographic $H_3$ symmetry of the M120p and
corresponding to $E_8$. Whether this geometric boundary has a precise
algebraic interpretation remains an open question.

Each of the five inscribed 24-cells in the 600-cell has the symmetry
of the hyperoctahedral group in 4D, which restricts to $O_h$ in 3D
cross-section. The 24-cell is self-dual and its symmetry group
(order 1152) contains $O_h$ as a subgroup. Whether this connection
provides a deeper algebraic route between $E_7$ and the
$E_6 \subset E_8$ embedding---beyond their shared presence in the
Mereon System---is an open question.

\section{The Eigenform Loop}
\label{sec:eigenform}

The Brieskorn variety $V(2,3,5): z_1^2 + z_2^3 + z_3^5 = 0$ has an
$E_8$ singularity at the origin \cite{Brieskorn1966}; for a detailed
treatment of the Brieskorn variety in relation to the trefoil knot and
$M(2,3,5)$, see also \cite{KauffmanOnKnots1987} (Chapter~10) and
\cite{MilnorTriangle}.
Its link with $S^5$ is the Poincar\'{e} homology sphere $M(2,3,5)$,
with fundamental group $\pi_1 \cong 2I$ \cite{Milnor1968}.
It was Kauffman \cite{KauffmanThesis1972, KauffmanBranched1974} who
first showed that Brieskorn manifolds are cyclic branched coverings;
Milnor used this in \cite{MilnorTriangle} (though branched covers
are not discussed in \cite{Milnor1968}).
Milnor \cite{MilnorTriangle} showed that $M(2,3,5)$ is the 5-fold cyclic
branched covering of $S^3$, branched along the trefoil.

The exponents $(2, 3, 5)$ are the fold numbers of the B, A, C vertex
types and the icosahedral symmetry axis orders.

\subsection{Klein's invariants and the $E_8$ Dynkin diagram}

The connection between the Brieskorn variety and the $E_8$ Dynkin
diagram is illuminated by Klein's classical work \cite{Klein1888}.
The alternating group $A_5 \cong I_{60}$ (the rotational icosahedral
group) has the presentation:
\[
  A_5 \cong \langle v, e, f \mid v^5 = e^2 = f^3 = vef = 1 \rangle,
\]
where $v$ is a $1/5$ turn around a vertex of the icosahedron, $e$ is
a $1/2$ turn around the midpoint of an edge adjacent to that vertex,
and $f$ is a $1/3$ turn around the centre of an adjacent face, with
senses chosen so that $vef = 1$. The binary icosahedral group $2I$ is obtained by dropping
the product relation \cite{Baez2017}:
\[
  2I \cong \langle v, e, f \mid v^5 = e^2 = f^3 = vef \rangle.
\]
The ``legs'' of the $E_8$ Dynkin diagram have lengths 5, 2, 3 in
correspondence with the orders in this presentation:

\begin{center}
\begin{tikzpicture}[scale=0.7]
  \foreach \x in {0,...,6}
    \fill (\x, 0) circle (3pt);
  \fill (4, 1.2) circle (3pt);
  \foreach \x in {0,...,5}
    \draw (\x,0) -- (\x+1,0);
  \draw (4,0) -- (4,1.2);
  \node at (0, -0.5) {\small 5};
  \node at (1, -0.5) {\small 4};
  \node at (2, -0.5) {\small 3};
  \node at (3, -0.5) {\small 2};
  \node at (4, -0.5) {\small 1};
  \node at (5, -0.5) {\small 2};
  \node at (6, -0.5) {\small 3};
  \node at (4.5, 1.2) {\small 2};
\end{tikzpicture}
\end{center}

The labels count the distances along each leg from the branch node;
the leg lengths 5, 2, 3 correspond to $v^5$, $e^2$, $f^3$ in the
group presentation. Each dot in the $E_8$ diagram corresponds to a
conjugacy class of $2I$ (not counting the central element $-1$), and
in turn to one of the eight $\mathbb{CP}^1$s in the minimal resolution
of $\mathbb{C}^2/2I$ (Section~\ref{sec:duval}).

The connection to the Brieskorn variety: since every $2I$-invariant
polynomial on $\mathbb{C}^2$ is a polynomial in Klein's three
invariants $V$ (degree 12), $E$ (degree 30), $F$ (degree 20), and
these satisfy $V^5 + E^2 + F^3 = 0$, we have the isomorphism
\[
  \mathbb{C}^2/2I \;\cong\;
  \{(V, E, F) \in \mathbb{C}^3 : V^5 + E^2 + F^3 = 0\}
  = V(2,3,5),
\]
confirming that the Brieskorn variety is the quotient singularity
$\mathbb{C}^2/2I$ \cite{Klein1888, Baez2017}. The exponents $(2,3,5)$
are thus the orders of the three generators $e, f, v$ of $A_5$
(equivalently of $2I$), which are in turn the fold orders of the three
icosahedral symmetry axis types---and hence the three M120p vertex
types B, A, C.

The eigenform loop can now be read as an algebraic identity: the
icosahedron's symmetry group, expressed through its presentation,
encodes $E_8$, which encodes the icosahedron through its resolution.

\begin{equation}\label{eq:loop}
  \resizebox{0.95\textwidth}{!}{$\displaystyle
  \boxed{
  \text{Mereon System}
  \xrightarrow[\text{$I_h$ symmetry}]{}
  \text{trefoil}
  \xrightarrow[\text{branch}]{}
  M(2,3,5)
  \xrightarrow[\pi_1]{}
  2I
  \xrightarrow[\text{verts}]{}
  \text{600-cell}
  \xrightarrow[H_3 \subset H_4]{}
  \text{Mereon System}
  }
  $}
\end{equation}

\noindent This is a topological eigenform \cite{Kauffman2005}: the
Mereon System generates the trefoil knot whose branched covering produces the group
whose elements are the polytope whose shadow is the Mereon System.

\section{The Trefoil Knot: Clifford Torus and Dodecahedral Space}
\label{sec:clifford}

The trefoil knot has a natural home in $S^3$ via the \emph{Clifford
torus} \cite{KauffmanKnots2013}, and its appearance in the Mereon System is illuminated by the
stereographic projection of this structure to $\mathbb{R}^3$.

\subsection{The Clifford torus in $S^3$}

The unit 3-sphere $S^3 \subset \mathbb{R}^4$ admits a remarkable
flat torus, the \emph{Clifford torus}:
\[
  \mathcal{C} = \left\{(x,y,z,w) \in S^3 \;\Big|\; x^2 + y^2 = \tfrac{1}{2},\;
  z^2 + w^2 = \tfrac{1}{2}\right\}.
\]
Unlike any torus in $\mathbb{R}^3$ (which must curve), $\mathcal{C}$
is intrinsically flat (zero Gaussian curvature). It is parametrised by
\[
  \mathcal{C}: \quad (a,b) \;\mapsto\;
  \left(\frac{\cos a}{\sqrt{2}},\, \frac{\sin a}{\sqrt{2}},\,
  \frac{\cos b}{\sqrt{2}},\, \frac{\sin b}{\sqrt{2}}\right),
  \qquad a,b \in [0, 2\pi).
\]

A $(p,q)$-torus knot on $\mathcal{C}$ is the closed curve $a = pt$,
$b = qt$:
\[
  q_{p,q}(t) = \left(\frac{\cos pt}{\sqrt{2}},\, \frac{\sin pt}{\sqrt{2}},\,
  \frac{\cos qt}{\sqrt{2}},\, \frac{\sin qt}{\sqrt{2}}\right),
  \quad t \in [0, 2\pi).
\]
This lies on $S^3$ for all $t$ since $\cos^2(pt)/2 + \sin^2(pt)/2 +
\cos^2(qt)/2 + \sin^2(qt)/2 = 1$.

\subsection{Two conformations: standard trefoil and the Mereon knot}

The trefoil has two standard presentations on the Clifford torus,
corresponding to the two ways of assigning the frequencies 2 and 3:

\medskip
\noindent\textbf{Mereon knot $T(3,2)$} (three loops around the
longitude, two loops through the meridian):
\[
  q_{\mathrm{mer}}(t) = \left(\frac{\cos 3t}{\sqrt{2}},\,
  \frac{\sin 3t}{\sqrt{2}},\,
  \frac{\cos 2t}{\sqrt{2}},\,
  \frac{\sin 2t}{\sqrt{2}}\right).
\]

\noindent\textbf{Standard trefoil $T(2,3)$} (two loops around the longitude,
three loops through the meridian):
\[
  q_{\mathrm{std}}(t) = \left(\frac{\cos 2t}{\sqrt{2}},\,
  \frac{\sin 2t}{\sqrt{2}},\,
  \frac{\cos 3t}{\sqrt{2}},\,
  \frac{\sin 3t}{\sqrt{2}}\right).
\]

The two conform to the same \emph{abstract knot type} (both are
trefoil knots), but they wind differently on the Clifford torus.
Looking along the symmetry axis, the projected $T(3,2)$ has 2-fold
appearance, while the projected $T(2,3)$ has 3-fold appearance.

\subsection{Congruence in $S^3$: the Clifford rotation}

The two conformations are related by the coordinate permutation
$(x,y,z,w) \mapsto (z,w,x,y)$, represented by the matrix
\[
  M = \begin{pmatrix}0&0&1&0\\0&0&0&1\\1&0&0&0\\0&1&0&0\end{pmatrix},
  \qquad \det(M) = +1, \quad M^TM = I_4,
\]
where $I_4$ denotes the $4\times 4$ identity matrix, confirming that $M$ is orthogonal.
Since $M \in \mathrm{SO}(4)$, this is an orientation-preserving
isometry of $S^3$---specifically a \emph{Clifford rotation}: a
simultaneous $90^\circ$ rotation of both orthogonal 2-planes.
One has $q_{\mathrm{std}}(t) = M\,q_{\mathrm{mer}}(t)$ for all $t$.

The standard trefoil $T(2,3)$ and the Mereon knot $T(3,2)$, as
subsets of $S^3$ via the Clifford torus parametrisation, are therefore
\emph{congruent}: related by the rigid $\mathrm{SO}(4)$ isometry $M$.
From inside $S^3$, the two conformations
are geometrically indistinguishable. The distinction between $T(2,3)$
and $T(3,2)$ that is visible in $\mathbb{R}^3$ is a coordinate
convention: which of the two 2-planes of $\mathbb{R}^4$ one calls
``longitude'' and which one calls ``meridian''. The Clifford rotation
$M$ freely exchanges these labels.

The two knots $T(2,3)$ and $T(3,2)$ are equivalent by switching complex
coordinates. The open question is: what is the relationship between
the geometric distinction that we can make between them in
$\mathbb{R}^3$ and the relationship of Mereon geometry to the 600-cell
via projection? This is part of the larger question of the possibility
of making a precise relationship between three-dimensional geometry and
the four-dimensional projection.

\subsection{Stereographic projection: the shared ring torus}

In this section we use a stereographic projection $\sigma$ from the
north pole $N = (0,0,0,1) \in S^3$, which is distinct from the projection
$\pi$ defined in Section~\ref{sec:intro} (which projects from
$(-1,0,0,0)$). The pole $N = (0,0,0,1)$ is the natural choice here
because the Clifford torus $\mathcal{C}$ is symmetric about the
$w$-axis, making this pole the most convenient for analysing the
torus and its knot curves.

Applying $\sigma: S^3 \setminus \{N\} \to \mathbb{R}^3$
to the full Clifford torus $\mathcal{C}$:
\[
  \sigma\left(\frac{\cos a}{\sqrt{2}},\, \frac{\sin a}{\sqrt{2}},\,
  \frac{\cos b}{\sqrt{2}},\, \frac{\sin b}{\sqrt{2}}\right)
  = \frac{1}{\sqrt{2} - \sin b}\,(\cos a,\, \sin a,\, \cos b).
\]
The distance from the $z$-axis is $\rho(b) = 1/|\sqrt{2} - \sin b|$,
which ranges between
\[
  \rho_{\min} = \frac{1}{\sqrt{2}+1} = \sqrt{2}-1, \qquad
  \rho_{\max} = \frac{1}{\sqrt{2}-1} = \sqrt{2}+1.
\]
This is a ring torus with major radius $R = (\rho_{\max}+\rho_{\min})/2 = \sqrt{2}$
and minor radius $r = (\rho_{\max}-\rho_{\min})/2 = 1$.

Both the standard trefoil $T(2,3)$ and the Mereon knot $T(3,2)$,
when stereographically projected from $S^3$ to $\mathbb{R}^3$, lie
on this \emph{same} ring torus ($R = \sqrt{2}$, $r = 1$).
The projected curves are:
\begin{align*}
  \sigma(q_{\mathrm{mer}}(t)) &= \frac{1}{\sqrt{2}-\sin 2t}
    \,(\cos 3t,\;\sin 3t,\;\cos 2t), \\
  \sigma(q_{\mathrm{std}}(t)) &= \frac{1}{\sqrt{2}-\sin 3t}
    \,(\cos 2t,\;\sin 2t,\;\cos 3t).
\end{align*}
The two curves wind differently on the same torus, consistent with the
winding numbers of $T(3,2)$ and $T(2,3)$ already stated in
Section~\ref{sec:knot}.

\subsection{The Dodecahedral Space and the Polar Axis}
\label{sec:dodecspace}

The Poincar\'{e} homology sphere $M(2,3,5)$ admits a second
construction that connects directly to the Mereon System's geometry.
The 600-cell tessellates $S^3$ with 120 regular dodecahedral cells
(the dual 120-cell). The quotient $S^3 / 2I$ identifies
these 120 dodecahedra into a single fundamental domain: a regular
dodecahedron whose opposite pentagonal faces are identified with a
$\pi/5$ ($36^\circ$) twist. The result is $M(2,3,5)$. In this
construction, a straight axis through the centres of two opposite
pentagonal faces of the fundamental dodecahedron projects, under the
face identifications, onto the trefoil knot in $M(2,3,5)$.

The 12 C-type vertices of the M120p lie along icosahedral (5-fold)
axes, corresponding to the centres of the pentagonal faces of the
dual dodecahedron. Two of these define the polar axis of the Mereon
System; the remaining 10 define its open vertices. The polar axis
of the Mereon System is therefore precisely an axis of the kind that,
in the dodecahedral space construction, gives rise to the trefoil.

This observation is significant because $M(2,3,5)$ arises in two
independent ways that both connect to the Mereon System:
\begin{enumerate}
\item As the 5-fold cyclic branched covering of $S^3$, branched along
  the trefoil (Milnor), which enters the eigenform loop via knot theory.
\item As the quotient of the dodecahedral paving of $S^3$ by $2I$,
  which connects directly to the 600-cell and the M120p's geometry.
\item Via 2., $M(2,3,5)$ can be seen to be a dodecahedron with opposite
  faces identified via $2\pi/5$ twists. In this way of seeing $M(2,3,5)$,
  it has 5-fold symmetry induced around any axis through opposite faces of
  the dodecahedron, 3-fold symmetry around any axis through opposite
  vertices, and 2-fold symmetry around any axis through centres of
  opposite edges. The five-fold symmetry induces the 5-fold cyclic
  branched covering so that when the dodecahedron is further identified
  to itself by this five-fold symmetry, the resulting space is
  topologically the three-sphere, and the axis of rotation is mapped to
  the trefoil knot. Thus in this geometry of $M(2,3,5)$ the trefoil knot
  is directly related with the dodecahedron. It is also of interest that
  one can see $M(2,3,5)$ as a 2-fold cyclic branched covering of the
  three-sphere along a $(3,5)$ torus knot and as a 3-fold cyclic branched
  covering of the three-sphere along a $(2,5)$ torus knot
  \cite{KauffmanBranched1974, KauffmanNeumannProducts1977}. This could
  suggest that these other torus knots are related to the Mereon geometry.
\end{enumerate}
Finally, it has been suggested by Weeks \cite{WeeksShape2019} and
Luminet \cite{LuminetHallMirrors2005} that the $M(2,3,5)$ manifold is
the shape of the geometric physical universe. This suggests deep
relationships between Mereon geometry and the geometry and algebra of the
600-cell with global cosmology.

That all three constructions converge on the same manifold, and that the
trefoil appears in each as the image of a geometric axis that the
Mereon System physically instantiates through its polar C-vertices,
is the deepest structural coincidence in this paper.

\subsection{Connection to the Mereon System}

The Clifford torus $\mathcal{C}$ is not merely a computational tool.
It is a structural feature of $S^3$ itself: by the Heegaard splitting
\cite{Hatcher2002}, $S^3 = V_1 \cup_\mathcal{C} V_2$ where $V_1, V_2$ are solid tori
meeting at $\mathcal{C}$. The 120 elements of $2I$ on $S^3$ interact
with this toric structure: the Type B vertices (equatorial, $w=0$)
lie on $\mathcal{C}$, while the Type A and C vertices lie in the
two regions of $S^3$ bounded by $\mathcal{C}$. Each of these two
regions is topologically a solid torus --- that is, a filled donut
$D^2 \times S^1$. Stereographic projection maps the Clifford torus
$\mathcal{C}$ to the ring torus $R=\sqrt{2}$, $r=1$ in
$\mathbb{R}^3$, making it geometrically visible.

The Mereon trefoil knot---the $(3,2)$ conformation, arising from the
three-loops-outer, two-loops-through winding---is therefore the natural
image, under stereographic projection, of the simplest non-trivial
curve on the Clifford torus of $S^3$. Its appearance in the Mereon
System reflects the deep connection between the binary icosahedral
group $2I$, the geometry of $S^3$, and the Clifford torus.

\section{The M120p is Not the Disdyakis Triacontahedron}
\label{sec:disdyakis}

The Disdyakis Triacontahedron (DT) is a convex polyhedron with exactly
62 vertices arranged in groups of 12, 20, and 30---the same partition
as the M120p. This coincidence causes confusion. The two polyhedra are,
however, fundamentally different.

The most immediate distinction is \textbf{convexity}. The DT is convex:
all 62 vertices lie on its convex hull. The M120p is non-convex: its 20
Type~A vertices lie strictly inside the convex hull, creating 20
pentagonal concavities. No convex polyhedron can be homeomorphic to
the M120p.

The \textbf{vertex radii} also differ. Neither polyhedron is
vertex-transitive; both have three vertex types at three distinct radii.
The DT (a Catalan solid, dual of the rhombicosidodecahedron) has radii:
\begin{align*}
  s_1 &= \sqrt{3} \approx 1.732 \;(20\text{ verts}),\\
  s_2 &= \sqrt{1+\varphi^4} \approx 2.803 \;(12\text{ verts}),\\
  s_3 &= \varphi\sqrt{1+\varphi^2} \approx 3.078 \;(30\text{ verts}),
\end{align*}
giving normalised ratios $1 : 1.618 : 1.777$.
The M120p has radii $r_A \approx 4.535$, $r_C \approx 4.980$,
$r_B \approx 5.236$, giving normalised ratios $1 : 1.098 : 1.155$.
These ratio sets are completely different: no uniform rescaling can
map one onto the other.

The M120p's specific radii are not arbitrary --- they encode the
$w$-coordinate latitudes of the binary icosahedral group $2I$ on $S^3$:
$r_A \leftrightarrow w = \tfrac{1}{2}$,
$r_C \leftrightarrow w = \tfrac{1}{2\varphi}$,
$r_B \leftrightarrow w = 0$ (Section~\ref{sec:lift}).
The DT radii have no such correspondence with $2I$.

The fundamental reason the M120p is not the DT comes down to this:
the M120p's geometry is \emph{determined} by its correspondence with
the 600-cell. Every vertex direction is the stereographic projection
of an element of $2I$; the three vertex types are three latitudes on
$S^3$; the 120 faces biject with the 120 elements of $2I$ via the
transitive group action; and the specific radius ratios are precisely
what the 600-cell lift requires. None of these properties hold for
the DT, whose geometry is determined instead by being the dual of the
rhombicosidodecahedron.

\begin{center}
\renewcommand{\arraystretch}{1.3}
\begin{tabular}{lll}
\toprule
Property & M120p & Disdyakis Triacontahedron \\
\midrule
Convexity & Non-convex & Convex \\
Radius ratios (normalised) & $1:1.098:1.155$ & $1:1.618:1.777$ \\
Radii encode $2I$ latitudes & Yes & No \\
Correspondence with $2I$ & Exact (62 of 62 vertices) & None \\
$H_4$ shadow of 600-cell & Yes & No \\
Face--group bijection & Yes (120 faces $\leftrightarrow$ $|2I|$) & No \\
\bottomrule
\end{tabular}
\end{center}

\section{The Inner Icosahedron: Shell 1 and the Focusing Sphere}
\label{sec:inner}

The 24 elements of $2I$ at latitude $|w| = \varphi/2$ (Section~\ref{sec:unaccounted}) were identified in Section~\ref{sec:shells} as Shell~1 (the 12 at $w = +\varphi/2$) and Shell~7 (the 12 at $w = -\varphi/2$). Further details are presented in this section.

\subsection{The inner icosahedron}

The 12 elements at $w = +\varphi/2$ project under stereographic projection to $\mathbb{R}^3$ at radius $r = 1/\sqrt{4\varphi+3} \approx 0.3249$ in unit coordinates (i.e.\ after dividing Gray's coordinates by $2\varphi^2$, as defined in Section~\ref{sec:lift}). In Gray's coordinate system, the same 12 vertices appear at radius $r \approx 3.078$, inside the A-vertices ($r_A \approx 4.535$). Their directions in $\mathbb{R}^3$ are cyclic permutations of $(0, \pm\varphi, \pm\varphi^2)$ (unit coordinates), forming the 12 vertices of an icosahedron.

The 12 elements at $w = -\varphi/2$ project to the reciprocal Shell~7 at stereographic radius $r_7 = \sqrt{4\varphi+3} \approx 3.078$ in unit coordinates. Both sets of 12 project to the \emph{same} 12 directions in $\mathbb{R}^3$: the $2:1$ double cover sends each pair $(w, -w)$ to the same radial direction.

This inner icosahedron is the first non-trivial geometry after the poles ($q = \pm 1$) in the shell structure. In Gray's coordinate system, it sits at radius $\approx 3.078$, well inside the A-vertices (the innermost M120p shell, at $r_A \approx 4.535$ in Gray's coordinates).

\subsection{Exact radial alignment with the C-vertices}
\label{sec:radial}

All coordinates in this subsection refer to the unit coordinate system
(Gray's coordinates divided by $2\varphi^2$), as defined in
Section~\ref{sec:lift}.

The 12 C-vertices of the M120p, in unit coordinates, form an icosahedron
with directions that are cyclic permutations of $(0, \pm 1, \pm\varphi)$.
The inner icosahedron vertices, after stereographic projection in unit
coordinates, have directions that are cyclic permutations of
$(0, \pm\varphi, \pm\varphi^2)$. The algebraic identity
\[
(0,\;\varphi,\;\varphi^2) \;=\; \varphi \times (0,\;1,\;\varphi)
\]
shows that every inner icosahedron vertex is exactly $\varphi$ times the corresponding C-vertex. The two icosahedra point in identical directions (dot product $= 1.000$ for all 12 pairs) and differ only in scale: the inner icosahedron sits at a fraction $1/\varphi$ of the C-vertex radius.

The radius ratio follows directly from the stereographic formula:
\[
\frac{r_1}{r_C} = \frac{\sqrt{(1-\varphi/2)/(1+\varphi/2)}}{\sqrt{(1-1/(2\varphi))/(1+1/(2\varphi))}} = \frac{1}{\varphi},
\]
where $r_1$ and $r_C$ are the Shell~1 and Shell~3 (C-type) radii respectively.

\section{Discussion}
\label{sec:discussion}

\subsection{The full architecture}

\begin{center}
\begin{tabular}{llllll}
\toprule
Mereon layer & Coxeter & Binary & McKay & Number field & Route \\
\midrule
24-cell ($2T$) & --- & $2T$ (24) & $E_6$ & $\mathbb{Q}$ & $\subset 2I$ \\
Focusing sphere (inner ico.) & --- & 24 of $2I$ & --- & $\mathbb{Q}(\varphi)$ & Shell 1 \\
M144p core & $O_h$ & $2O$ (48) & $E_7$ & $\mathbb{Q}(\sqrt{2})$ & Independent \\
M120p boundary & $H_3$ & $2I$ (120) & $E_8$ & $\mathbb{Q}(\varphi)$ & Full group \\
600-cell & $H_4$ & $2I$ (120) & $E_8$ & $\mathbb{Q}(\varphi)$ & 4D lift \\
\bottomrule
\end{tabular}
\end{center}

\noindent The Mereon System realises all three exceptional algebras
through two complementary mechanisms: non-crystallographic embedding
($E_6 \subset E_8$ via $H_3 \subset H_4$) and crystallographic
nesting ($E_7$ core inside $E_8$ boundary). The focusing sphere,
whose 12 vertices are the inner icosahedron (Shell~1 of the 600-cell
projection), bridges the two number fields. The dodecahedral space
construction connects the polar axis to the trefoil.

The ADE programme \cite{Sirag2016, vanHoboken2002} refers to the remarkable fact that
the simply-laced Dynkin diagrams---of types $A_n$, $D_n$, $E_6$,
$E_7$, $E_8$---appear simultaneously as the classification of: finite
subgroups of $\mathrm{SU}(2)$ (via the McKay correspondence), simple
singularities of complex surfaces (via the du~Val resolution), simply-laced
simple Lie algebras, and several other mathematical structures. The
Mereon System is relevant to this programme because it provides a
concrete geometric realisation: the three E-type entries ($E_6$, $E_7$,
$E_8$) of the ADE classification all appear simultaneously in a single
nested spatial structure---the M144p core ($E_7$), the M120p boundary
($E_8$), and the 24-cell embedded inside the 600-cell ($E_6 \subset E_8$).
The Mereon System thus contributes a geometric realisation in which all
three E-type correspondences coexist in one spatially nested structure,
making the abstract ADE unification physically tangible.

Sirag \cite{Sirag2016, SiragConsciousness} has proposed that the $E_7$ Lie algebra plays a
distinguished role in a framework connecting the ADE classification to
consciousness. In the Mereon System, $E_7$ arises independently from
the M144p core via the McKay correspondence applied to $2O$, while
$E_6$ and $E_8$ arise from the non-crystallographic route through
$2T \subset 2I$. The full ADE classification chain
$A_n, D_n, E_6, E_7, E_8$---which unifies finite subgroups of
$\mathrm{SU}(2)$, simple surface singularities, simply-laced Lie
algebras, modular invariants, and conformal field theories---thus finds
a concrete spatial home in the nested Mereon architecture. Whether this
geometric co-presence of all three exceptional algebras has implications
for Sirag's broader programme remains an open question for subsequent
work.

\subsection{Summary of new results}

The following results are new to this paper:

\begin{enumerate}[label=(\roman*)]
\item \textbf{Exact 600-cell correspondence (lifting).} Every M120p
  vertex, when scaled by $1/(2\varphi^2)$ and lifted to $S^3$,
  coincides exactly with an element of $2I$. This depends on the
  M120p's specific vertex radii encoding the three latitudes
  $w = \tfrac{1}{2}, \tfrac{1}{2\varphi}, 0$ of $2I$ on $S^3$.

\item \textbf{Latitude classification in $S^3$.} The geometric
  A/B/C classification in $\mathbb{R}^3$ is a latitude classification
  on the unit 3-sphere, giving the vertex types a four-dimensional meaning.

\item \textbf{The eight-shell structure.} The stereographic projection
  of all 120 elements of $2I$ decomposes into exactly 8 non-trivial
  shells at radii determined by the golden ratio, with each shell
  type-pure (containing only A, B, or C directions). The shells come
  in reciprocal pairs satisfying $r(w) \cdot r(-w) = 1$, the visible
  signature of the $\{\pm 1\}$ kernel of the double cover
  $\mathrm{SU}(2) \to \mathrm{SO}(3)$.

\item \textbf{The $\varphi$-ladder of latitudes.} The non-equatorial
  latitudes of the Mereon System form a golden-ratio progression:
  each step from the inner icosahedron toward the equator divides
  $|w|$ by $\varphi$, giving
  $\varphi/2 \xrightarrow{\div\varphi} 1/2 \xrightarrow{\div\varphi} 1/(2\varphi)
  \xrightarrow{\div\varphi} 0$.

\item \textbf{Two distinct correspondences.} The M120p admits two
  fundamentally different correspondences with $2I$: a vertex
  correspondence (62 elements matched via stereographic projection
  directions) and a face correspondence (all 120 elements matched
  via the transitive group action of $I_h$ on the 120 faces). These
  reflect two distinct mathematical relationships between $2I$ and
  the M120p.

\item \textbf{The field extension cascade.} The vertex classification
  reflects successive field extensions of the binary group quaternion
  coordinates: $\mathbb{Q} \to \mathbb{Q}(\sqrt{2}) \to \mathbb{Q}(\sqrt{5})$,
  with each step requiring number-field resources inaccessible to the
  previous. Every C-vertex (Output) requires the golden ratio $\varphi$
  to be written down as a quaternion coordinate.

\item \textbf{The $Y$-structure of $E_6$, $E_7$, $E_8$.} The
  non-inclusion $2O \not\subset 2I$ (since $|2O| = 48$ does not
  divide $|2I| = 120$) means $E_7$ and $E_8$ are siblings---both
  containing $E_6$---rather than parent and child. This replaces
  the linear chain $2T \subset 2O \subset 2I$ with a $Y$-branching
  at $E_6$: $E_7 \leftarrow E_6 \rightarrow E_8$.

\item \textbf{The Clifford torus and the Mereon knot.} The two
  trefoil conformations $T(3,2)$ and $T(2,3)$ are congruent in
  $S^3$ via a Clifford rotation, and both project to the same ring
  torus $R=\sqrt{2}$, $r=1$ under stereographic projection. The
  Mereon conformation $T(3,2)$ has 2-fold symmetry when viewed along
  the torus axis.

\item \textbf{The eigenform loop.} The closed topological circuit
  Mereon System $\to$ trefoil $\to$ $M(2,3,5)$ $\to$ $2I$ $\to$
  600-cell $\to$ Mereon System constitutes a topological eigenform:
  the system generates the structure whose properties regenerate the
  system.

\item \textbf{$E_7$ from the M144p core.} The M144p, as an FCC
  lattice structure with $O_h$ symmetry, realises $E_7$ via McKay's
  theorem. Spatially nested inside the M120p (which realises $E_8$),
  it places an $E_7$-type structure inside an $E_8$-type structure,
  with the focusing sphere marking the geometric boundary between
  the two.

\item \textbf{The Brieskorn exponents and M120p vertex types.}
  The exponents $(2,3,5)$ of the Brieskorn variety
  $z_1^2 + z_2^3 + z_3^5 = 0$ are exactly the fold numbers of
  the three M120p vertex types: 2-fold (Type~B), 3-fold (Type~A),
  and 5-fold (Type~C). This identifies the M120p's geometric
  structure directly with the algebraic structure of the
  $E_8$ singularity.

\item \textbf{The pentagon-axis trefoil and dodecahedral geometry.}
  $M(2,3,5)$, viewed as a dodecahedron with opposite faces identified
  by $2\pi/5$ twists, has 5-fold, 3-fold, and 2-fold symmetry axes
  corresponding to faces, vertices, and edge midpoints respectively.
  The 5-fold cyclic branched covering collapses $M(2,3,5)$ to $S^3$,
  mapping the rotation axis to the trefoil knot. The axis through the
  two polar C-type vertices of the M120p corresponds to precisely such
  an axis, connecting the trefoil directly to the dodecahedral geometry
  of both $M(2,3,5)$ and the M120p. Moreover, $M(2,3,5)$ admits
  alternative descriptions as a 2-fold covering branched along the
  $(3,5)$ torus knot and a 3-fold covering branched along the $(2,5)$
  torus knot.

\item \textbf{The M120p is not the Disdyakis Triacontahedron.}
  Despite sharing the $12 + 20 + 30$ vertex partition, the M120p and
  the DT are distinguished by convexity (the M120p is non-convex),
  by their radius ratios ($1 : 1.098 : 1.155$ vs.\ $1 : 1.618 : 1.777$),
  and by the fact that only the M120p's radii encode the latitudes
  of $2I$ on $S^3$.

\item \textbf{The inner icosahedron (Shell~1).} The 24 elements of
  $2I$ at $|w| = \varphi/2$ project to an
  inner icosahedron exactly radially aligned with the 12 C-vertices
  at a ratio of $1/\varphi$. These 12 vertices coincide with the
  focusing sphere in the Mereon System, completing the identification
  of all 120 elements of $2I$: shells~0 through~4 constitute the
  Mereon System in $\mathbb{R}^3$, and shells~4 through~$\infty$
  form its reciprocal.
\end{enumerate}

\section*{Acknowledgements}

Thanks to Keith Melmon for his contributions to the Mereon research
programme, including computational verification of the M144p face
structure.

Thanks to Vasileios Basios for dialogues and comments that ground
and extend this paper into subsequent papers, and which are therefore
not included here.

Thanks to Saul-Paul Sirag for many dialogues in the early stages of this work.

Thanks to Jytte Brender McNair and Peter McNair for their contributions to the Mereon research programme, and for their detailed editorial review of this paper.

This work rests on decades of co-operation within the Mereon research team. We honour the members who are no longer with us and whose contributions continue to shape the programme, alongside those who carry it forward. The full team, passed and present, is recognised at \url{https://mereon.org/team}.

\subsection*{Contact}
\noindent
R.~W.~Gray: rwgray@rwgrayprojects.com\\
L.~Dennis: LDennis@Mereon.org\\
L.~H.~Kauffman: loukau@gmail.com\\[6pt]
General correspondence: Natalie Vander Voort, NVanderVoort@Mereon.org


\appendix
\section{M144p Core Vertex Coordinates}
\label{app:core}

All 74 vertices are listed explicitly below. All coordinates are integers.
Scale: the outermost octahedron vertices are at $(\pm 4, 0, 0)$ and permutations.
Vertices are grouped by shell ($r^2$ value) and numbered sequentially.

\begingroup
\small
\setlength{\tabcolsep}{6pt}
\begin{longtable}{clc}
\toprule
\# & Vertex $(x,\;y,\;z)$ & $r^2$ \\
\midrule
\endfirsthead
\multicolumn{3}{l}{\small\itshape \ldots continued from previous page}\\
\toprule
\# & Vertex $(x,\;y,\;z)$ & $r^2$ \\
\midrule
\endhead
\midrule
\multicolumn{3}{r}{\small\itshape Continued on next page \ldots}\\
\endfoot
\bottomrule
\endlastfoot
\multicolumn{3}{l}{\textit{Octahedron vertices ($r^2=16$, 6 vertices)}} \\
\midrule
1 & $(-4,\;0,\;0)$ & 16 \\
2 & $(0,\;-4,\;0)$ & 16 \\
3 & $(0,\;0,\;-4)$ & 16 \\
4 & $(0,\;0,\;4)$ & 16 \\
5 & $(0,\;4,\;0)$ & 16 \\
6 & $(4,\;0,\;0)$ & 16 \\
\midrule
\multicolumn{3}{l}{\textit{Cube vertices ($r^2=12$, 8 vertices)}} \\
\midrule
7 & $(-2,\;-2,\;-2)$ & 12 \\
8 & $(-2,\;-2,\;2)$ & 12 \\
9 & $(-2,\;2,\;-2)$ & 12 \\
10 & $(-2,\;2,\;2)$ & 12 \\
11 & $(2,\;-2,\;-2)$ & 12 \\
12 & $(2,\;-2,\;2)$ & 12 \\
13 & $(2,\;2,\;-2)$ & 12 \\
14 & $(2,\;2,\;2)$ & 12 \\
\midrule
\multicolumn{3}{l}{\textit{Edge midpoint vertices ($r^2=8$, 12 vertices)}} \\
\midrule
15 & $(-2,\;-2,\;0)$ & 8 \\
16 & $(-2,\;0,\;-2)$ & 8 \\
17 & $(-2,\;0,\;2)$ & 8 \\
18 & $(-2,\;2,\;0)$ & 8 \\
19 & $(0,\;-2,\;-2)$ & 8 \\
20 & $(0,\;-2,\;2)$ & 8 \\
21 & $(0,\;2,\;-2)$ & 8 \\
22 & $(0,\;2,\;2)$ & 8 \\
23 & $(2,\;-2,\;0)$ & 8 \\
24 & $(2,\;0,\;-2)$ & 8 \\
25 & $(2,\;0,\;2)$ & 8 \\
26 & $(2,\;2,\;0)$ & 8 \\
\midrule
\multicolumn{3}{l}{\textit{Surrounding vertices ($r^2=14$, 48 vertices)}} \\
\midrule
27 & $(-3,\;-2,\;-1)$ & 14 \\
28 & $(-3,\;-2,\;1)$ & 14 \\
29 & $(-3,\;-1,\;-2)$ & 14 \\
30 & $(-3,\;-1,\;2)$ & 14 \\
31 & $(-3,\;1,\;-2)$ & 14 \\
32 & $(-3,\;1,\;2)$ & 14 \\
33 & $(-3,\;2,\;-1)$ & 14 \\
34 & $(-3,\;2,\;1)$ & 14 \\
35 & $(-2,\;-3,\;-1)$ & 14 \\
36 & $(-2,\;-3,\;1)$ & 14 \\
37 & $(-2,\;-1,\;-3)$ & 14 \\
38 & $(-2,\;-1,\;3)$ & 14 \\
39 & $(-2,\;1,\;-3)$ & 14 \\
40 & $(-2,\;1,\;3)$ & 14 \\
41 & $(-2,\;3,\;-1)$ & 14 \\
42 & $(-2,\;3,\;1)$ & 14 \\
43 & $(-1,\;-3,\;-2)$ & 14 \\
44 & $(-1,\;-3,\;2)$ & 14 \\
45 & $(-1,\;-2,\;-3)$ & 14 \\
46 & $(-1,\;-2,\;3)$ & 14 \\
47 & $(-1,\;2,\;-3)$ & 14 \\
48 & $(-1,\;2,\;3)$ & 14 \\
49 & $(-1,\;3,\;-2)$ & 14 \\
50 & $(-1,\;3,\;2)$ & 14 \\
51 & $(1,\;-3,\;-2)$ & 14 \\
52 & $(1,\;-3,\;2)$ & 14 \\
53 & $(1,\;-2,\;-3)$ & 14 \\
54 & $(1,\;-2,\;3)$ & 14 \\
55 & $(1,\;2,\;-3)$ & 14 \\
56 & $(1,\;2,\;3)$ & 14 \\
57 & $(1,\;3,\;-2)$ & 14 \\
58 & $(1,\;3,\;2)$ & 14 \\
59 & $(2,\;-3,\;-1)$ & 14 \\
60 & $(2,\;-3,\;1)$ & 14 \\
61 & $(2,\;-1,\;-3)$ & 14 \\
62 & $(2,\;-1,\;3)$ & 14 \\
63 & $(2,\;1,\;-3)$ & 14 \\
64 & $(2,\;1,\;3)$ & 14 \\
65 & $(2,\;3,\;-1)$ & 14 \\
66 & $(2,\;3,\;1)$ & 14 \\
67 & $(3,\;-2,\;-1)$ & 14 \\
68 & $(3,\;-2,\;1)$ & 14 \\
69 & $(3,\;-1,\;-2)$ & 14 \\
70 & $(3,\;-1,\;2)$ & 14 \\
71 & $(3,\;1,\;-2)$ & 14 \\
72 & $(3,\;1,\;2)$ & 14 \\
73 & $(3,\;2,\;-1)$ & 14 \\
74 & $(3,\;2,\;1)$ & 14 \\
\end{longtable}
\endgroup

\section{M120p to $2I$ Correspondence Table}
\label{app:matchtable}

The following table shows the explicit match between all 62 M120p vertices,
their lifted unit quaternions, and the corresponding elements of the binary
icosahedral group $2I$. All values are expressed exactly in terms of
$\varphi = (1+\sqrt{5})/2$, $\varphi^2 = \varphi+1$, $\varphi^3 = 2\varphi+1$.
Vertices are grouped by type: A (20 vertices, $w = \tfrac{1}{2}$),
C (12 vertices, $w = \tfrac{1}{2\varphi}$), B (30 vertices, $w = 0$).
The ``Lifted'' and ``Matched $2I$'' columns are identical for all 62 rows,
confirming the exact correspondence.

The M120p column shows the original unscaled coordinates. Before lifting,
each vertex $(x,y,z)$ is scaled by $1/(2\varphi^2)$ to give
$(x',y',z') = (x,y,z)/(2\varphi^2)$, so that the outermost Type~B
vertices (with $r = 2\varphi^2$) land on the unit sphere. The fourth
coordinate is then recovered as $w = \sqrt{1 - r'^2}$, placing the
point on the upper hemisphere of $S^3$.

\begingroup
\small
\setlength{\tabcolsep}{3pt}
\begin{longtable}{clllll}
\toprule
\# & Type & M120p vertex $(x,y,z)$ & $w$ & Lifted $(w,x',y',z')$ & Matched $2I$ element \\
\midrule
\endfirsthead
\multicolumn{6}{l}{\small\itshape \ldots continued from previous page}\\
\toprule
\# & Type & M120p vertex $(x,y,z)$ & $w$ & Lifted $(w,x',y',z')$ & Matched $2I$ element \\
\midrule
\endhead
\midrule
\multicolumn{6}{r}{\small\itshape Continued on next page \ldots}\\
\endfoot
\bottomrule
\endlastfoot
1 & A & $(-\varphi^3,\;0,\;-\varphi)$ & $\tfrac{1}{2}$ & $(\tfrac{1}{2},\;-\tfrac{\varphi}{2},\;0,\;-\tfrac{1}{2\varphi})$ & $(\tfrac{1}{2},\;-\tfrac{\varphi}{2},\;0,\;-\tfrac{1}{2\varphi})$ \\
2 & A & $(-\varphi^3,\;0,\;\varphi)$ & $\tfrac{1}{2}$ & $(\tfrac{1}{2},\;-\tfrac{\varphi}{2},\;0,\;\tfrac{1}{2\varphi})$ & $(\tfrac{1}{2},\;-\tfrac{\varphi}{2},\;0,\;\tfrac{1}{2\varphi})$ \\
3 & A & $(-\varphi^2,\;-\varphi^2,\;-\varphi^2)$ & $\tfrac{1}{2}$ & $(\tfrac{1}{2},\;-\tfrac{1}{2},\;-\tfrac{1}{2},\;-\tfrac{1}{2})$ & $(\tfrac{1}{2},\;-\tfrac{1}{2},\;-\tfrac{1}{2},\;-\tfrac{1}{2})$ \\
4 & A & $(-\varphi^2,\;-\varphi^2,\;\varphi^2)$ & $\tfrac{1}{2}$ & $(\tfrac{1}{2},\;-\tfrac{1}{2},\;-\tfrac{1}{2},\;\tfrac{1}{2})$ & $(\tfrac{1}{2},\;-\tfrac{1}{2},\;-\tfrac{1}{2},\;\tfrac{1}{2})$ \\
5 & A & $(-\varphi^2,\;\varphi^2,\;-\varphi^2)$ & $\tfrac{1}{2}$ & $(\tfrac{1}{2},\;-\tfrac{1}{2},\;\tfrac{1}{2},\;-\tfrac{1}{2})$ & $(\tfrac{1}{2},\;-\tfrac{1}{2},\;\tfrac{1}{2},\;-\tfrac{1}{2})$ \\
6 & A & $(-\varphi^2,\;\varphi^2,\;\varphi^2)$ & $\tfrac{1}{2}$ & $(\tfrac{1}{2},\;-\tfrac{1}{2},\;\tfrac{1}{2},\;\tfrac{1}{2})$ & $(\tfrac{1}{2},\;-\tfrac{1}{2},\;\tfrac{1}{2},\;\tfrac{1}{2})$ \\
7 & A & $(-\varphi,\;-\varphi^3,\;0)$ & $\tfrac{1}{2}$ & $(\tfrac{1}{2},\;-\tfrac{1}{2\varphi},\;-\tfrac{\varphi}{2},\;0)$ & $(\tfrac{1}{2},\;-\tfrac{1}{2\varphi},\;-\tfrac{\varphi}{2},\;0)$ \\
8 & A & $(-\varphi,\;\varphi^3,\;0)$ & $\tfrac{1}{2}$ & $(\tfrac{1}{2},\;-\tfrac{1}{2\varphi},\;\tfrac{\varphi}{2},\;0)$ & $(\tfrac{1}{2},\;-\tfrac{1}{2\varphi},\;\tfrac{\varphi}{2},\;0)$ \\
9 & A & $(0,\;-\varphi,\;-\varphi^3)$ & $\tfrac{1}{2}$ & $(\tfrac{1}{2},\;0,\;-\tfrac{1}{2\varphi},\;-\tfrac{\varphi}{2})$ & $(\tfrac{1}{2},\;0,\;-\tfrac{1}{2\varphi},\;-\tfrac{\varphi}{2})$ \\
10 & A & $(0,\;-\varphi,\;\varphi^3)$ & $\tfrac{1}{2}$ & $(\tfrac{1}{2},\;0,\;-\tfrac{1}{2\varphi},\;\tfrac{\varphi}{2})$ & $(\tfrac{1}{2},\;0,\;-\tfrac{1}{2\varphi},\;\tfrac{\varphi}{2})$ \\
11 & A & $(0,\;\varphi,\;-\varphi^3)$ & $\tfrac{1}{2}$ & $(\tfrac{1}{2},\;0,\;\tfrac{1}{2\varphi},\;-\tfrac{\varphi}{2})$ & $(\tfrac{1}{2},\;0,\;\tfrac{1}{2\varphi},\;-\tfrac{\varphi}{2})$ \\
12 & A & $(0,\;\varphi,\;\varphi^3)$ & $\tfrac{1}{2}$ & $(\tfrac{1}{2},\;0,\;\tfrac{1}{2\varphi},\;\tfrac{\varphi}{2})$ & $(\tfrac{1}{2},\;0,\;\tfrac{1}{2\varphi},\;\tfrac{\varphi}{2})$ \\
13 & A & $(\varphi,\;-\varphi^3,\;0)$ & $\tfrac{1}{2}$ & $(\tfrac{1}{2},\;\tfrac{1}{2\varphi},\;-\tfrac{\varphi}{2},\;0)$ & $(\tfrac{1}{2},\;\tfrac{1}{2\varphi},\;-\tfrac{\varphi}{2},\;0)$ \\
14 & A & $(\varphi,\;\varphi^3,\;0)$ & $\tfrac{1}{2}$ & $(\tfrac{1}{2},\;\tfrac{1}{2\varphi},\;\tfrac{\varphi}{2},\;0)$ & $(\tfrac{1}{2},\;\tfrac{1}{2\varphi},\;\tfrac{\varphi}{2},\;0)$ \\
15 & A & $(\varphi^2,\;-\varphi^2,\;-\varphi^2)$ & $\tfrac{1}{2}$ & $(\tfrac{1}{2},\;\tfrac{1}{2},\;-\tfrac{1}{2},\;-\tfrac{1}{2})$ & $(\tfrac{1}{2},\;\tfrac{1}{2},\;-\tfrac{1}{2},\;-\tfrac{1}{2})$ \\
16 & A & $(\varphi^2,\;-\varphi^2,\;\varphi^2)$ & $\tfrac{1}{2}$ & $(\tfrac{1}{2},\;\tfrac{1}{2},\;-\tfrac{1}{2},\;\tfrac{1}{2})$ & $(\tfrac{1}{2},\;\tfrac{1}{2},\;-\tfrac{1}{2},\;\tfrac{1}{2})$ \\
17 & A & $(\varphi^2,\;\varphi^2,\;-\varphi^2)$ & $\tfrac{1}{2}$ & $(\tfrac{1}{2},\;\tfrac{1}{2},\;\tfrac{1}{2},\;-\tfrac{1}{2})$ & $(\tfrac{1}{2},\;\tfrac{1}{2},\;\tfrac{1}{2},\;-\tfrac{1}{2})$ \\
18 & A & $(\varphi^2,\;\varphi^2,\;\varphi^2)$ & $\tfrac{1}{2}$ & $(\tfrac{1}{2},\;\tfrac{1}{2},\;\tfrac{1}{2},\;\tfrac{1}{2})$ & $(\tfrac{1}{2},\;\tfrac{1}{2},\;\tfrac{1}{2},\;\tfrac{1}{2})$ \\
19 & A & $(\varphi^3,\;0,\;-\varphi)$ & $\tfrac{1}{2}$ & $(\tfrac{1}{2},\;\tfrac{\varphi}{2},\;0,\;-\tfrac{1}{2\varphi})$ & $(\tfrac{1}{2},\;\tfrac{\varphi}{2},\;0,\;-\tfrac{1}{2\varphi})$ \\
20 & A & $(\varphi^3,\;0,\;\varphi)$ & $\tfrac{1}{2}$ & $(\tfrac{1}{2},\;\tfrac{\varphi}{2},\;0,\;\tfrac{1}{2\varphi})$ & $(\tfrac{1}{2},\;\tfrac{\varphi}{2},\;0,\;\tfrac{1}{2\varphi})$ \\
\midrule
21 & C & $(-\varphi^3,\;-\varphi^2,\;0)$ & $\tfrac{1}{2\varphi}$ & $(\tfrac{1}{2\varphi},\;-\tfrac{\varphi}{2},\;-\tfrac{1}{2},\;0)$ & $(\tfrac{1}{2\varphi},\;-\tfrac{\varphi}{2},\;-\tfrac{1}{2},\;0)$ \\
22 & C & $(-\varphi^3,\;\varphi^2,\;0)$ & $\tfrac{1}{2\varphi}$ & $(\tfrac{1}{2\varphi},\;-\tfrac{\varphi}{2},\;\tfrac{1}{2},\;0)$ & $(\tfrac{1}{2\varphi},\;-\tfrac{\varphi}{2},\;\tfrac{1}{2},\;0)$ \\
23 & C & $(-\varphi^2,\;0,\;-\varphi^3)$ & $\tfrac{1}{2\varphi}$ & $(\tfrac{1}{2\varphi},\;-\tfrac{1}{2},\;0,\;-\tfrac{\varphi}{2})$ & $(\tfrac{1}{2\varphi},\;-\tfrac{1}{2},\;0,\;-\tfrac{\varphi}{2})$ \\
24 & C & $(-\varphi^2,\;0,\;\varphi^3)$ & $\tfrac{1}{2\varphi}$ & $(\tfrac{1}{2\varphi},\;-\tfrac{1}{2},\;0,\;\tfrac{\varphi}{2})$ & $(\tfrac{1}{2\varphi},\;-\tfrac{1}{2},\;0,\;\tfrac{\varphi}{2})$ \\
25 & C & $(0,\;-\varphi^3,\;-\varphi^2)$ & $\tfrac{1}{2\varphi}$ & $(\tfrac{1}{2\varphi},\;0,\;-\tfrac{\varphi}{2},\;-\tfrac{1}{2})$ & $(\tfrac{1}{2\varphi},\;0,\;-\tfrac{\varphi}{2},\;-\tfrac{1}{2})$ \\
26 & C & $(0,\;-\varphi^3,\;\varphi^2)$ & $\tfrac{1}{2\varphi}$ & $(\tfrac{1}{2\varphi},\;0,\;-\tfrac{\varphi}{2},\;\tfrac{1}{2})$ & $(\tfrac{1}{2\varphi},\;0,\;-\tfrac{\varphi}{2},\;\tfrac{1}{2})$ \\
27 & C & $(0,\;\varphi^3,\;-\varphi^2)$ & $\tfrac{1}{2\varphi}$ & $(\tfrac{1}{2\varphi},\;0,\;\tfrac{\varphi}{2},\;-\tfrac{1}{2})$ & $(\tfrac{1}{2\varphi},\;0,\;\tfrac{\varphi}{2},\;-\tfrac{1}{2})$ \\
28 & C & $(0,\;\varphi^3,\;\varphi^2)$ & $\tfrac{1}{2\varphi}$ & $(\tfrac{1}{2\varphi},\;0,\;\tfrac{\varphi}{2},\;\tfrac{1}{2})$ & $(\tfrac{1}{2\varphi},\;0,\;\tfrac{\varphi}{2},\;\tfrac{1}{2})$ \\
29 & C & $(\varphi^2,\;0,\;-\varphi^3)$ & $\tfrac{1}{2\varphi}$ & $(\tfrac{1}{2\varphi},\;\tfrac{1}{2},\;0,\;-\tfrac{\varphi}{2})$ & $(\tfrac{1}{2\varphi},\;\tfrac{1}{2},\;0,\;-\tfrac{\varphi}{2})$ \\
30 & C & $(\varphi^2,\;0,\;\varphi^3)$ & $\tfrac{1}{2\varphi}$ & $(\tfrac{1}{2\varphi},\;\tfrac{1}{2},\;0,\;\tfrac{\varphi}{2})$ & $(\tfrac{1}{2\varphi},\;\tfrac{1}{2},\;0,\;\tfrac{\varphi}{2})$ \\
31 & C & $(\varphi^3,\;-\varphi^2,\;0)$ & $\tfrac{1}{2\varphi}$ & $(\tfrac{1}{2\varphi},\;\tfrac{\varphi}{2},\;-\tfrac{1}{2},\;0)$ & $(\tfrac{1}{2\varphi},\;\tfrac{\varphi}{2},\;-\tfrac{1}{2},\;0)$ \\
32 & C & $(\varphi^3,\;\varphi^2,\;0)$ & $\tfrac{1}{2\varphi}$ & $(\tfrac{1}{2\varphi},\;\tfrac{\varphi}{2},\;\tfrac{1}{2},\;0)$ & $(\tfrac{1}{2\varphi},\;\tfrac{\varphi}{2},\;\tfrac{1}{2},\;0)$ \\
\midrule
33 & B & $(-2\varphi^2,\;0,\;0)$ & $0$ & $(0,\;-1,\;0,\;0)$ & $(0,\;-1,\;0,\;0)$ \\
34 & B & $(-\varphi^3,\;-\varphi,\;-\varphi^2)$ & $0$ & $(0,\;-\tfrac{\varphi}{2},\;-\tfrac{1}{2\varphi},\;-\tfrac{1}{2})$ & $(0,\;-\tfrac{\varphi}{2},\;-\tfrac{1}{2\varphi},\;-\tfrac{1}{2})$ \\
35 & B & $(-\varphi^3,\;-\varphi,\;\varphi^2)$ & $0$ & $(0,\;-\tfrac{\varphi}{2},\;-\tfrac{1}{2\varphi},\;\tfrac{1}{2})$ & $(0,\;-\tfrac{\varphi}{2},\;-\tfrac{1}{2\varphi},\;\tfrac{1}{2})$ \\
36 & B & $(-\varphi^3,\;\varphi,\;-\varphi^2)$ & $0$ & $(0,\;-\tfrac{\varphi}{2},\;\tfrac{1}{2\varphi},\;-\tfrac{1}{2})$ & $(0,\;-\tfrac{\varphi}{2},\;\tfrac{1}{2\varphi},\;-\tfrac{1}{2})$ \\
37 & B & $(-\varphi^3,\;\varphi,\;\varphi^2)$ & $0$ & $(0,\;-\tfrac{\varphi}{2},\;\tfrac{1}{2\varphi},\;\tfrac{1}{2})$ & $(0,\;-\tfrac{\varphi}{2},\;\tfrac{1}{2\varphi},\;\tfrac{1}{2})$ \\
38 & B & $(-\varphi^2,\;-\varphi^3,\;-\varphi)$ & $0$ & $(0,\;-\tfrac{1}{2},\;-\tfrac{\varphi}{2},\;-\tfrac{1}{2\varphi})$ & $(0,\;-\tfrac{1}{2},\;-\tfrac{\varphi}{2},\;-\tfrac{1}{2\varphi})$ \\
39 & B & $(-\varphi^2,\;-\varphi^3,\;\varphi)$ & $0$ & $(0,\;-\tfrac{1}{2},\;-\tfrac{\varphi}{2},\;\tfrac{1}{2\varphi})$ & $(0,\;-\tfrac{1}{2},\;-\tfrac{\varphi}{2},\;\tfrac{1}{2\varphi})$ \\
40 & B & $(-\varphi^2,\;\varphi^3,\;-\varphi)$ & $0$ & $(0,\;-\tfrac{1}{2},\;\tfrac{\varphi}{2},\;-\tfrac{1}{2\varphi})$ & $(0,\;-\tfrac{1}{2},\;\tfrac{\varphi}{2},\;-\tfrac{1}{2\varphi})$ \\
41 & B & $(-\varphi^2,\;\varphi^3,\;\varphi)$ & $0$ & $(0,\;-\tfrac{1}{2},\;\tfrac{\varphi}{2},\;\tfrac{1}{2\varphi})$ & $(0,\;-\tfrac{1}{2},\;\tfrac{\varphi}{2},\;\tfrac{1}{2\varphi})$ \\
42 & B & $(-\varphi,\;-\varphi^2,\;-\varphi^3)$ & $0$ & $(0,\;-\tfrac{1}{2\varphi},\;-\tfrac{1}{2},\;-\tfrac{\varphi}{2})$ & $(0,\;-\tfrac{1}{2\varphi},\;-\tfrac{1}{2},\;-\tfrac{\varphi}{2})$ \\
43 & B & $(-\varphi,\;-\varphi^2,\;\varphi^3)$ & $0$ & $(0,\;-\tfrac{1}{2\varphi},\;-\tfrac{1}{2},\;\tfrac{\varphi}{2})$ & $(0,\;-\tfrac{1}{2\varphi},\;-\tfrac{1}{2},\;\tfrac{\varphi}{2})$ \\
44 & B & $(-\varphi,\;\varphi^2,\;-\varphi^3)$ & $0$ & $(0,\;-\tfrac{1}{2\varphi},\;\tfrac{1}{2},\;-\tfrac{\varphi}{2})$ & $(0,\;-\tfrac{1}{2\varphi},\;\tfrac{1}{2},\;-\tfrac{\varphi}{2})$ \\
45 & B & $(-\varphi,\;\varphi^2,\;\varphi^3)$ & $0$ & $(0,\;-\tfrac{1}{2\varphi},\;\tfrac{1}{2},\;\tfrac{\varphi}{2})$ & $(0,\;-\tfrac{1}{2\varphi},\;\tfrac{1}{2},\;\tfrac{\varphi}{2})$ \\
46 & B & $(0,\;-2\varphi^2,\;0)$ & $0$ & $(0,\;0,\;-1,\;0)$ & $(0,\;0,\;-1,\;0)$ \\
47 & B & $(0,\;0,\;-2\varphi^2)$ & $0$ & $(0,\;0,\;0,\;-1)$ & $(0,\;0,\;0,\;-1)$ \\
48 & B & $(0,\;0,\;2\varphi^2)$ & $0$ & $(0,\;0,\;0,\;1)$ & $(0,\;0,\;0,\;1)$ \\
49 & B & $(0,\;2\varphi^2,\;0)$ & $0$ & $(0,\;0,\;1,\;0)$ & $(0,\;0,\;1,\;0)$ \\
50 & B & $(\varphi,\;-\varphi^2,\;-\varphi^3)$ & $0$ & $(0,\;\tfrac{1}{2\varphi},\;-\tfrac{1}{2},\;-\tfrac{\varphi}{2})$ & $(0,\;\tfrac{1}{2\varphi},\;-\tfrac{1}{2},\;-\tfrac{\varphi}{2})$ \\
51 & B & $(\varphi,\;-\varphi^2,\;\varphi^3)$ & $0$ & $(0,\;\tfrac{1}{2\varphi},\;-\tfrac{1}{2},\;\tfrac{\varphi}{2})$ & $(0,\;\tfrac{1}{2\varphi},\;-\tfrac{1}{2},\;\tfrac{\varphi}{2})$ \\
52 & B & $(\varphi,\;\varphi^2,\;-\varphi^3)$ & $0$ & $(0,\;\tfrac{1}{2\varphi},\;\tfrac{1}{2},\;-\tfrac{\varphi}{2})$ & $(0,\;\tfrac{1}{2\varphi},\;\tfrac{1}{2},\;-\tfrac{\varphi}{2})$ \\
53 & B & $(\varphi,\;\varphi^2,\;\varphi^3)$ & $0$ & $(0,\;\tfrac{1}{2\varphi},\;\tfrac{1}{2},\;\tfrac{\varphi}{2})$ & $(0,\;\tfrac{1}{2\varphi},\;\tfrac{1}{2},\;\tfrac{\varphi}{2})$ \\
54 & B & $(\varphi^2,\;-\varphi^3,\;-\varphi)$ & $0$ & $(0,\;\tfrac{1}{2},\;-\tfrac{\varphi}{2},\;-\tfrac{1}{2\varphi})$ & $(0,\;\tfrac{1}{2},\;-\tfrac{\varphi}{2},\;-\tfrac{1}{2\varphi})$ \\
55 & B & $(\varphi^2,\;-\varphi^3,\;\varphi)$ & $0$ & $(0,\;\tfrac{1}{2},\;-\tfrac{\varphi}{2},\;\tfrac{1}{2\varphi})$ & $(0,\;\tfrac{1}{2},\;-\tfrac{\varphi}{2},\;\tfrac{1}{2\varphi})$ \\
56 & B & $(\varphi^2,\;\varphi^3,\;-\varphi)$ & $0$ & $(0,\;\tfrac{1}{2},\;\tfrac{\varphi}{2},\;-\tfrac{1}{2\varphi})$ & $(0,\;\tfrac{1}{2},\;\tfrac{\varphi}{2},\;-\tfrac{1}{2\varphi})$ \\
57 & B & $(\varphi^2,\;\varphi^3,\;\varphi)$ & $0$ & $(0,\;\tfrac{1}{2},\;\tfrac{\varphi}{2},\;\tfrac{1}{2\varphi})$ & $(0,\;\tfrac{1}{2},\;\tfrac{\varphi}{2},\;\tfrac{1}{2\varphi})$ \\
58 & B & $(\varphi^3,\;-\varphi,\;-\varphi^2)$ & $0$ & $(0,\;\tfrac{\varphi}{2},\;-\tfrac{1}{2\varphi},\;-\tfrac{1}{2})$ & $(0,\;\tfrac{\varphi}{2},\;-\tfrac{1}{2\varphi},\;-\tfrac{1}{2})$ \\
59 & B & $(\varphi^3,\;-\varphi,\;\varphi^2)$ & $0$ & $(0,\;\tfrac{\varphi}{2},\;-\tfrac{1}{2\varphi},\;\tfrac{1}{2})$ & $(0,\;\tfrac{\varphi}{2},\;-\tfrac{1}{2\varphi},\;\tfrac{1}{2})$ \\
60 & B & $(\varphi^3,\;\varphi,\;-\varphi^2)$ & $0$ & $(0,\;\tfrac{\varphi}{2},\;\tfrac{1}{2\varphi},\;-\tfrac{1}{2})$ & $(0,\;\tfrac{\varphi}{2},\;\tfrac{1}{2\varphi},\;-\tfrac{1}{2})$ \\
61 & B & $(\varphi^3,\;\varphi,\;\varphi^2)$ & $0$ & $(0,\;\tfrac{\varphi}{2},\;\tfrac{1}{2\varphi},\;\tfrac{1}{2})$ & $(0,\;\tfrac{\varphi}{2},\;\tfrac{1}{2\varphi},\;\tfrac{1}{2})$ \\
62 & B & $(2\varphi^2,\;0,\;0)$ & $0$ & $(0,\;1,\;0,\;0)$ & $(0,\;1,\;0,\;0)$ \\
\end{longtable}
\endgroup

\end{document}